\documentclass[11pt]{amsart}

\usepackage{graphicx}
\usepackage{amsfonts}
\usepackage{mathrsfs}
\usepackage[colorlinks=true, pdfstartview=FitV, linkcolor=blue, citecolor=blue, urlcolor=blue]{hyperref}
\usepackage{color}
\usepackage{graphicx,epstopdf,verbatim}
\usepackage[cal=boondoxo]{mathalfa}

\allowdisplaybreaks

\usepackage{graphicx}
\usepackage{amsmath}
\usepackage{amsthm}
\usepackage{xcolor}
\usepackage{bbm}

 \usepackage{amssymb,amscd,stmaryrd}

\newtheorem{lem}{Lemma}
\newtheorem{prop}{Proposition}
\newtheorem{thm}{Theorem}
\newtheorem{coro}{Corollary}
\newtheorem{remark}{Remark}

\newcommand{\R}{\mathbb{R}}

\begin{document}

\title{Two-weight estimates for the square function and $t$-Haar multipliers}
\author{Daewon Chung, Jean Carlo Moraes, Mar\'ia Cristina Pereyra, Brett Wick}
\date{September 2024}

\dedicatory{Dedicated to Jill Pipher for everything she has done as a mathematician and mentor, building community,  and as a tireless advocate for mathematics. You are an inspiration to all of us!}

\begin{abstract}
We present necessary and sufficient conditions on triples of weights $(u,v,w)$ for the boundedness of the dyadic weighted square function $S_w$ from $L^2(u)$ into $L^2(v)$. We use this characterization to obtain  necessary and sufficient conditions for the boundedness of the $t$-Haar multipliers from $L^2(u)$ into $L^2(v)$ in terms of boundedness of the dyadic weighted square function.
\end{abstract} 


\maketitle

\tableofcontents

\section{Introduction}
The \emph{dyadic weighted square function} is defined for locally integrable functions $f$  formally~as
\[ S^d_wf(x)= \bigg (\sum_{I\in\mathcal{D}} |\langle f\rangle_I^w-\langle f\rangle_{\widetilde{I}}^w|^2 \mathbbm{1}_I(x)\bigg )^{1/2},\]
where $w$ is a weight, $\mathcal{D}$ denotes the dyadic intervals, $\widetilde{I}$ is the parent of $I\in\mathcal{D}$, $\mathbbm{1}_I$ denotes the characteristic function of the set $I$,  and $\langle f\rangle_I^w=\frac{\int_I f(x)w(x)\,dx }{\int_I w(x)\, dx}$  denotes the weighted integral average of $f$ with respect to  $w$. 

When $w=1$, this corresponds to the  dyadic square function $S^d$, known to be bounded on $L^p(\mathbb{R})$ for all $1<p<\infty$, and bounded on weighted $L^p(u)$ if and only if $u\in A_p^d$ (the Muckenhoupt dyadic $A_p$ class of weights), see for example \cite{BL79,Buc93}. Necessary and sufficient conditions  for the two-weight boundedness of $S^d$ were found in \cite{NTV99}. The following quantitative estimates were obtained, using Bellman functions,  in \cite{HTV00, Wit00} for the linear upper bound,  and in \cite{Wil89, PP02} for the lower bound,
\[ c[u]_{A_2^d}^{-1/2}\|f\|_{L^2(u)} \leq \|S^df\|_{L^2(u)} \leq C[u]_{A_2^d}\|f\|_{L^2(u)}.\]
In the matrix valued setting, the $A_2$ conjecture in the matrix case (the one-weight linear upper bound) was proven in \cite{HPV19}, but so far  the square-root lower bound  is off by a logarithm factor  \cite{BPW16}.  In a non-homogeneous setting it was shown, in \cite{DIPTV19}, that the square-root lower bound has to be replaced by a linear bound. In \cite{Tr23}, using  sparse domination, a sharp  linear mixed $A_2$-$A_{\infty}$  two-weight upper bound estimate in the matrix case was found for the vector square function  in the non-homogeneous setting.

In this paper, we are concerned with quantitative  two-weight upper bounds for $S_w$ in the scalar setting. 
Given weights $u,v, w$ that are locally integrable and a.e. non-negative functions on $\mathbb{R}$. 
We present  sufficient conditions for the boundedness of the dyadic weighted square function $S_w$ from $L^2(u)$ into $L^2(v)$. That is, we provide conditions on the triple of weights $(u,v,w)$ so that there is a constant $C=C(u,v,w)>0$ such that for all $f\in L^2(u)$ we have
\[ \|S_w^df\|_{L^2(v)}\leq C \|f\|_{L^2(u)}.\]
 If, in addition, we assume that $w$ is a dyadic doubling weight, the sufficient conditions are also necessary. We present  necessary conditions (without assuming doubling) for the two weight boundedness of the dyadic weighted square function. These results recover previous  qualitative known results.  In order to be precise, we present the results in question that will be proved in Section~\ref{SquareFunction}. In this paper,  $I^+$ and $I^-$ denote the right and left children of the dyadic interval $I$, respectively, and $w(I)=\int_Iw(x)\,dx$.

\begin{thm}[The sufficiency]\label{thm:1} Let $u,v,w$ be weights, such that  $u^{-1}w^2,v\in L^1_{{\rm loc}}(\mathbb{R})$. If
\begin{itemize}
    \item[{\rm (i)}] there is a constant $C_1>0$ such that for all $J\in\mathcal{D}$ we have that 
    $$  \langle u^{-1}w \rangle_J^w \langle vw^{-1} \rangle_J^w \leq C_1;$$
    
     \item[{\rm (ii)}] there is a $C_2>0$ such that for all $J\in\mathcal{D}$ we have that,
$$\frac{1}{w(J)} \sum_{I \in \mathcal{D}(J) } \big| \Delta^w_I (u^{-1}w) \big|^2  \left [ \frac{w(I^-)^2v(I^+) + w(I^+)^2v(I^-)}{w(I)^2} \right ]  \leq C_2 \langle u^{-1}w\rangle_J^w,$$
where $\Delta_I^w(g)=\langle g\rangle_{I^+}^w-\langle g\rangle_{I^-}^w$.
\end{itemize}
Then, the weighted dyadic square function $S^d_w$
is bounded from $L^2(u)$ into $L^2(v)$. 
Furthermore,
\[\|S_w^d\|_{L^2(u)\to L^2(v)} \leq \sqrt{C_1}+2\sqrt{C_2}.\]
\end{thm}

We cannot recover conditions (i) and (ii) as necessary by requiring only that \( u, v, w \) are weights satisfying \( u^{-1}w^2, v \in L^1_{\mathrm{loc}}(\mathbb{R}) \). However, this becomes possible in some situations, for example if we additionally assume that \( w \) is a dyadic doubling weight. In order to state our necessary conditions, let us denote, for a given $J\in\mathcal{D}$, $J^*$ as its dyadic sibling, so that $J\cup J^*=\widetilde{J}$.

\begin{thm}[The necessity]\label{thm:2} Let $u,v,w$ be weights, such that  $u^{-1}w^2,v\in L^1_{{\rm loc}}(\mathbb{R})$. Suppose $S_w$ is bounded from $L^2(u)$ into $L^2(v)$ then   
\begin{itemize}
    \item[{\rm (i)}]  There is a constant $C>0$ such that for all $J\in\mathcal{D}$ we have that 
    $$  \left(\frac{w(J^*)}{w(\widetilde{J})}\right)^2 \langle u^{-1}w \rangle^w_J \langle vw^{-1} \rangle^w_J \leq C.$$
     \item[{\rm (ii)}]  There is a constant $C>0$ such that for all $J\in\mathcal{D}$ we have that
$$\frac{1}{w(J)} \sum_{I \in \mathcal{D}(J) } \big| \Delta^w_I (u^{-1}w) \big|^2  \left [ \frac{w(I^-)^2v(I^+) + w(I^+)^2v(I^-)}{w(I)^2} \right ]  \leq C \langle u^{-1}w\rangle_J^w.$$
\end{itemize}
\end{thm}
 


We say that a pair of weights $(u,v)\in A_2^d(w)$ if and only if 
\[[u,v]_{A_2^d(w)}:=\sup_{I \in \mathcal{D}} \langle u^{-1}\rangle_I^w \langle v\rangle_I^w< \infty.\]
Condition (i) in Theorem~\ref{thm:1}  can be rephrased as $(uw^{-1},vw^{-1})\in A_2^d(w)$. Condition (i) in Theorem~\ref{thm:2} leads us to introduce a  new class of weights, that we call the {\it restricted weighted dyadic $A_2$ class}, denoted $\tilde{A}_2^d(w)$. A pair of weights $(u,v)\in \tilde{A}_2^d(w)$ if and only if 
\[[u,v]_{\tilde{A}_2^d(w)}:=\sup_{I \in \mathcal{D}}  \left(\frac{w(J^*)}{w(\widetilde{J})}\right)^2\langle u^{-1}\rangle_I^w \langle v\rangle_I^w< \infty.\]
Condition (i) in Theorem~\ref{thm:2} can be rephrased as $(uw^{-1},vw^{-1}) \in \tilde{A}_2^d(w)$.

Note that 
 $A_2^d(w)\subset \tilde{A}_2^d(w)$, and when $w$ is dyadic doubling the classes coincide.

As an application, we obtain  necessary and sufficient conditions for the boundedness of the $t$-Haar multipliers from $L^2(u)$ into $L^2(v)$, but under some additional conditions on the triple of weights $(u,v,w)$.

The \emph{(variable) Haar multiplier associated to a weight $w$} is defined by
\[T_wf(x)=\sum_{I\in\mathcal{D}} \frac{w(x)}{\langle w\rangle_I}\langle f,h_I\rangle h_I,\]
where $\{h_I: I\in\mathcal{D}\}$ are the Haar functions, an orthogonal basis in $L^2(\R)$ and an unconditional basis in $L^p(\R )$, and $\langle f,g\rangle$ denotes the inner product in $L^2(\R)$.
These Haar multipliers
appeared in connection to the resolvent of the dyadic paraproduct in work by the third author \cite{Per94}. This connection was established thanks to a celebrated product formula studied by Jill Pipher and collaborators in \cite{FKP91}. We will say more about this connection in Section~\ref{sec:tHaar}. It was shown in \cite{Per94}, that for $1<p<\infty$, the operator $T_w$ is bounded on $L^p(\R)$ if and only if $S_w^d$ is bounded on $L^{p'}(\R )$ with $\frac1p+\frac{1}{p'}=1$, and either of these happen  if, and only if, $w$ belongs to the dyadic Reverse H\"older-$p$ class, $w\in RH^d_p$, namely, 
\begin{equation}\label{def:RHp}
[w]_{RH_p^d}:=\sup_{I\in \mathcal{D}} \frac{\langle w^p\rangle_I^{1/p}}{\langle w\rangle_I} <\infty.
\end{equation}
Note that when $u=v=1$, condition (i) in Theorem~\ref{thm:1} 
is precisely $w\in RH_2^d$. Although not obvious, in this case, condition (ii)   is also equivalent to $w\in RH_2^d$, and that will imply that condition (i) in Theorem~\ref{thm:2} also holds, as  $[1,1]_{\tilde{A}_2^d(w)} \leq [w]_{RH_2^d}=[1,1]_{A_2^d(w)}$.  So that in this case, we recover the result ``$S_w$ is bounded on $L^2(\R )$ if, and only if, $w\in RH_2^d$." 

Having recently  identified necessary and sufficient conditions for the two-weight boundedness of a larger class of Haar multipliers \cite{CHMPW24}, the signed $t$-Haar multipliers, we wondered if we could have necessary and sufficient two-weight boundedness conditions for the signed $t$-Haar multipliers involving boundedness of the dyadic weighted square function as the result of \cite{Per94} previously mentioned.

The \emph{signed $t$-Haar multipliers}, $T^t_{w,\sigma}$ are defined for a weight $w$, a real number $t\in\R$, and a sequence of signs $\sigma=\{\sigma_I:I\in\mathcal{D}\}$ where $\sigma_I=\pm1$, by
\begin{equation}\label{eq:t-HaarMultiplier}
T^t_{w,\sigma}f(x)=\sum_{I\in\mathcal{D}} \sigma_I \left (\frac{w(x)}{\langle w\rangle_I}\right )^t \langle f,h_I\rangle h_I.
\end{equation}
When $\sigma_I=1$ for all $I\in\mathcal{D}$ and $t=1$ this is $T_w$. When $w=1$ or equivalently $t=0$ then this is the martingale transform $T_{\sigma}$. We  have results for $T_{w,\sigma}$, the case $t=1$, that we describe in detail in Section~\ref{sec:tHaar}.\\

%
%
%
%

In Section~\ref{sec:prelim} we introduce the basic definitions of  dyadic intervals, weights, dyadic doubling weights, dyadic Muckenhoupt $A_p^d$, $RH_p^d$, $A_{\infty}^d$, and $RH_1^d$ classes of weights and their characterization by summation conditions. We also introduce  weighted Haar functions, weighted Carleson sequences, and Sawyer's estimates. In Section~\ref{sec:3weight}, we introduce three-weight Muckenhoupt like conditions and a larger restricted class of weights that appears naturally in our two-weight theorems for the weighted square function.
In Section~\ref{sec:Sw}, we prove Theorem~\ref{thm:1} and Theorem~\ref{thm:2}, the sufficient and necessary conditions for the two-weight boundedness of the weighted dyadic square function.
In Section~\ref{sec:previous}, we verify that previous known results for the weighted and unweighted dyadic square function are recovered. Finally, in Section~\ref{sec:tHaar}, we state and prove a theorem connecting the two-weight boundedness of the Haar multiplier $T_{w,\sigma}$ to the two-weight boundedness of the weighted square function $S_w$. We also describe how the Haar multiplier $T_w$ is connected to the resolvent of the dyadic paraproduct via the celebrated Fefferman-Kenig-Pipher product formula.\\

\noindent{\bf Acknowledgements.} The second author was partially supported by the Coordenação de Aperfeiçoamento de Pessoal de Nível Superior – Brasil (CAPES) – Edital 41/2017 (Chamada 001/2024 – PROPG/UFRGS).  The last author was partially supported by ARC DP 220100285 as well as National Science Foundation Grants DMS \# 2054863 and \# 2349868.

\section{Preliminaries}\label{sec:prelim}
In this section, we present the concepts of weights, dyadic intervals, dyadic \(A^d_p\) and \(RH^d_p\) weight classes for $1<p<\infty$, and the dyadic $A_{\infty}^d$ and $RH_1^d$ weight classes, the Haar basis and its weighted counterparts, as well as Carleson sequences. Additionally, we summarize the celebrated Fefferman-Kenig-Pipher estimates for $A_{\infty}^d$ weights, and Buckley's summation conditions for $A_p^d$, $RH_p^d$, and $RH_1^d$ weights, as well as a  lesser known characterization of $RH_p^d$, all needed in the subsequent discussion. We also state Beznosova's Little Lemma for Carleson sequences and a weighted Sawyer's estimate that will be used later in the paper.

\subsection{Dyadic Intervals and Weights}\label{sec:weights}

The family of \emph{dyadic intervals}, denoted $\mathcal{D}$, consists of  those intervals of the form $I=[k2^{-j},(k+1)2^{-j})$ for $k,j\in\mathbb{Z}$.  The dyadic intervals are nested, given $I,J\in\mathcal{D}$ then exactly one of the following holds: $I\cap J=\emptyset$, $I\subsetneq J$, $I=J$, or $J\subsetneq I$.
We denote $\mathcal{D}(J)$ the collection of dyadic subintervals of a dyadic  interval $J$, that is,  $\mathcal{D}(J)=\{I\in\mathcal{D}: I\subset J\}$ for all $J\in\mathcal{D}$. Each dyadic interval $I$ has exactly one \emph{parent}, denoted $\widetilde{I}$; one \emph{sibling}, denoted $I^*$; and two \emph{children}, denoted $I^+$ and $I^-$, the right and left halves respectively of $I$. These are all dyadic intervals, so that $I^+\cup I^- = I \subsetneq \widetilde{I}=I\cup I^*$. Necessarily, $|\widetilde{I}|=2|I|$ and $|I|=2|I^{\pm}|$.

A function  $\omega:\R\to \R$ is a \emph{weight} if it is an almost everywhere positive function that is locally integrable. 

We use the notation $\langle \omega\rangle_I$ to denote the integral average of $\omega$ on $I$, that is $\frac{1}{|I|}\int_I\omega(x)\,dx$, and $|I|$ denotes the length of the interval $I$. We will also use the notation $\omega(I)=\int_I \omega(x)\,dx$.

A weight $\omega$ is \emph{dyadic doubling} if, and only if, $D(\omega):=\sup_{I\in\mathcal{D}} \omega (\widetilde{I})/\omega(I) <\infty$, and  $D(\omega)$ is called the \emph{doubling constant of }$\omega$.

A weight $\omega$ belongs to the \emph{dyadic $A_p$ Muckenhoupt class} for $1<p<\infty$ if and only if 
\[ [\omega]_{A_p^d}:=\sup_{I\in\mathcal{D}} \langle \omega\rangle_I \langle \omega^{-\frac{1}{p-1}}\rangle_I^{p-1} <\infty.\]
The quantity $[\omega ]_{A_p^d}$ is called the $A_p^d$-characteristic of the weight $\omega$.

Note that $\omega\in A^d_p$ if and only if $\omega^{-1/(p-1)} \in A^d_{p'}$ where $\frac1p+\frac{1}{p'}=1$. In particular $\omega\in A_2^d$ if and only if $\omega^{-1}\in A_2^d$.

A weight $\omega$ belongs to the \emph{dyadic Reverse H\"older-$p$ class}, $RH_p^d$,  for $1<p<\infty$ if and only if \eqref{def:RHp} holds, in other words there is a constant $C>0$ such that for all $I\in\mathcal{D}$,
\[  \langle \omega^p\rangle_I^{1/p} \leq C \langle \omega \rangle_I.\]
The smallest constant $C$ is the \emph{$RH_p^d$-characteristic} of the weight $\omega$, denoted $[\omega]_{RH_p^d}$.

The following characterization of $A_p^d$ and $RH_p^d$ weights by summation conditions is due to Stephen Buckley.
\begin{prop}[\cite{Buc93}]\label{prop:Buckley} Let $1<p<\infty$ and $\omega$ a dyadic doubling weight. 

{\rm (i)}  The weight $\omega\in A_p^d$ if, and only if, there is a $C>0$ such that for all $J\in\mathcal{D}$
\[ \frac{1}{|J|} \sum_{I\in \mathcal{D}(J)} |I|  \left |\frac{\Delta_I \omega}{\langle \omega\rangle_I}\right |^2\langle \omega\rangle_I^{-\frac{1}{p-1}} \leq C \langle \omega\rangle_J^{-\frac{1}{p-1}}.\]

{\rm (ii)} The weight  $\omega$ is in $RH_p^d$ if, and only if, there is a $C>0$ such that for all $J\in\mathcal{D}$
\[ \frac{1}{|J|} \sum_{I\in \mathcal{D}(J)} |I|  \left |\frac{\Delta_I \omega}{\langle \omega\rangle_I}\right |^2\langle \omega\rangle_I^p \leq C \langle \omega\rangle_J^p.\]
Where we recall $\Delta_I\omega:= \langle \omega\rangle_{I^+}-\langle \omega\rangle_{I^-}$.
\end{prop}

Buckley's result is an outgrowth of the celebrated Fefferman-Kenig-Pipher characterization of dyadic $A_{\infty}$, that appeared in \cite{FKP91}. 

A weight $\omega$ belongs to $ A_{\infty}^d$ if, and only if, 
$$ [\omega]_{A_{\infty}^d} := \sup_{I\in\mathcal{D}} \langle \omega\rangle_I e^{-\langle \log \omega\rangle_I}<\infty.$$
It is well known  that $\omega\in A_{\infty}^d$ implies that $\omega $ is dyadic doubling, and that $A_{\infty}^d=\cup_{p>1} A_p^d$.

\begin{prop}[\cite{FKP91}, Theorem 3.22] \label{thm:FKP}
A weight $w\in A_{\infty}^d$ if, and only if,  there is a constant $C>0$ such that for every dyadic interval $J$,
\[ \frac{1}{|J|} \sum_{I\in \mathcal{D}(J)} \left |\frac{\Delta_Iw}{\langle w\rangle_I}\right |^2 |I| \leq C \log [w]_{A_{\infty}^d}.\]
\end{prop}

A weight $\omega\in RH_1^d$ if, and only if, 
\[ [\omega]_{RH_1^d}=\sup_{I\in\mathcal{D}} \left \langle \frac{\omega}{\langle \omega\rangle_I} \log \frac{\omega}{\langle \omega\rangle_I} \right \rangle_I <\infty.\]

\begin{prop}[\cite{Buc93}]\label{prop:RH1}
 The weight $\omega\in RH_1^d$ if, and only if, there is $C>0$ such that for all $J\in\mathcal{D}$
\[ \frac{1}{|J|} \sum_{I\in \mathcal{D}(J)} |I|  \left |\frac{\Delta_I \omega}{\langle \omega\rangle_I}\right |^2\langle \omega\rangle_I \leq C \langle \omega\rangle_J.\]
\end{prop} 
The class $RH_1^d$ coincides with $A_{\infty}^d$ when the weights are dyadic doubling. In general,
$A_{\infty}^d\subset RH_1^d$. In \cite[Theorem I.5]{BR14}, it is shown that  $[\omega]_{RH_1^d} \leq (\ln 16) \, [\omega]_{A_{\infty}^d}$, using Bellman functions.

 We will also need another characterization of $RH_p^d$ that appeared in \cite[Proposition 1]{CHMPW24}, which  is a consequence of Buckley's characterization of $RH_1^d$ (Proposition~\ref{prop:RH1}), and the fact that for $1<p<\infty$, $w\in RH_p^d$ if, and only if, $w^p\in RH_1^d$, that can be found in  \cite[Lemma 1]{KP99}.

\begin{prop}[\cite{CHMPW24}]\label{lem:CHMPW}
Let $1<p<\infty$. A weight $w\in RH_p^d$ if and only if there is a constant $C>0$ such that
\[ \frac{1}{|J|} \sum_{I\in\mathcal{D}(J)} |I|\,  \frac{|\Delta_I(w^p)|^2}{\langle w\rangle^p_I}  \leq C \langle w^p\rangle_J, \quad\quad \mbox{for all $J\in\mathcal{D}$.}\]
\end{prop}

\subsection{Haar Functions}\label{sec:Haar}
Given a weight $\omega$, and an interval $I$, the \emph{weighted Haar function associated to $I$} is defined as follows,
\[h_I^\omega(x)= \sqrt{\frac{\omega(I^-)}{\omega(I^+)\omega(I)}}\mathbbm{1}_{I^+}(x) - \sqrt{\frac{\omega(I^+)}{\omega(I^-)\omega(I)}}\mathbbm{1}_{I^-}(x).\]
Where the constants have been chosen to enforce $L^2(\omega)$-norm one and mean zero, that is: $\|h_I^\omega\|_{L^2(\omega)}=1$ and $\int_Ih_I^\omega(x)\omega(x)\,dx=0$. These properties, coupled with  the nested properties of the dyadic intervals, ensure that $\{h_I^{\omega}\}_{I\in\mathcal{D}}$ forms an orthonormal system in $L^2(\omega)$. In particular, Bessel's inequality holds, that is for $f\in L^2(\omega)$ then
\[ \sum_{I\in\mathcal{D}} |\langle f,h_I^{\omega}\rangle_{\omega}|^2 \leq \|f\|_{L^2(\omega)}^2,\]
where $\langle f,g\rangle_{\omega}=\int_{\R} f(x)\,g(x)\,\omega(x)\,dx$, the inner product in $L^2(\omega)$.

When $\omega=1$, we recover the standard Haar functions $h_I=|I|^{-1/2}(\mathbbm{1}_{I^+}-\mathbbm{1}_{I^-})$. The collection $\{h_I\}_{I\in\mathcal{D}}$ forms  an orthonormal basis in $L^2(\R )$, and an unconditional basis in $L^2(\omega)$ if, and only if, $w\in A_2^d$ \cite{GW74}.

We can always write $h_I^\omega$ as a linear combination of another weighted Haar function, say $h_I^{\nu}$ and the characteristic function $\mathbbm{1}_I$. More precisely,
\begin{equation}\label{eqn:weightedHaar-linear-combination}
h_I^\omega=\alpha_I^{\nu,\omega} h_I^{\nu} + \beta_I^{\nu,\omega} \mathbbm{1}_I,
\end{equation}
where 
\[ \alpha_I^{\nu,\omega}= \sqrt{\frac{\langle \nu \omega^{-1}\rangle_{I^+}^\omega\langle \nu \omega^{-1}\rangle_{I^-}^\omega}{\langle \nu \omega^{-1}\rangle_{I}^\omega}}, \quad \beta_I^{\nu,\omega} = \sqrt{\frac{\omega(I^+)\,\omega(I^-)}{\omega(I)}} \frac{1}{\omega(I)} \frac{\Delta_I^\omega(\nu \omega^{-1})}{\langle \nu \omega^{-1}\rangle_I^\omega},\]
with $\Delta_I^\omega(\nu )= \langle \nu \rangle_{I^+}^\omega -\langle \nu\rangle_{I^-}^\omega$, and $\langle \nu\rangle_I^\omega=\langle \nu \omega\rangle_I/\langle \omega\rangle_I$ is the integral average of $\nu$ with respect to $\omega$.

Note that when $\omega=1$ then $|\langle f\rangle^\omega_I-\langle f\rangle^\omega_{\widetilde{I}}|=\frac12 |\Delta^\omega_{\widetilde{I}}(f)|$, but for a general weight this is no longer true. What is always true is that 
\begin{align}
  |\langle f\rangle^\omega_I  - \langle f\rangle^\omega_{\widetilde{I}} |  
  &= \sqrt{\frac{\omega(I^*)}{\omega(\widetilde{I})\omega(I)}} \big| \langle f, h^\omega_{\widetilde{I}} \rangle_\omega \big| \, = \, \big|\langle f,h^\omega_{\widetilde{I}}\rangle_\omega h^\omega_{\widetilde{I}}(x_I)\big| \quad\quad \mbox{where $x_I\in I$}. \label{eq:Delta^w}
\end{align}

Another calculation shows that
\begin{equation}\label{eq:DeltaIw(f)}
|\langle f, h_I^\omega\rangle_\omega|^2 = \frac{\omega(I^+) \omega(I^-)}{\omega(I)}\big |\langle f\rangle_{I^+}^\omega-\langle f\rangle_{I^-}^\omega\big |^2 =  \frac{\omega(I^+) \omega(I^-)}{\omega(I)}\big |\Delta_I^\omega(f)\big |^2.
\end{equation}

\subsection{Weighted Carleson Sequences}\label{sec:Carleson}
Given a sequence of real numbers $\{\lambda_I\}_{I\in\mathcal{D}}$, we say that it is a \emph{Carleson sequence} if, and only if,  it satisfies the following packing or Carleson condition: there is a constant $B>0$ such that for all $J\in \mathcal{D}$ we have
\begin{equation}\label{Carleson-sequence}
\frac{1}{|J|} \sum_{I\in\mathcal{D}(J)} \lambda_I \leq B. 
\end{equation}

Given a weight $\omega$, we say that the sequence $\{\lambda_I\}_{I\in\mathcal{D}}$ is a \emph{$\omega$-Carleson sequence}, if, and only if, there is a constant $B>0$ such that for all $J\in \mathcal{D}$ we have
\begin{equation}\label{w-Carleson-sequence}
\frac{1}{|J|} \sum_{I\in\mathcal{D}(J)} \lambda_I \langle \omega \rangle_I\leq B \langle \omega \rangle_J.
\end{equation}
We are following here the definition of a weighted Carleson sequence as in \cite{Per94}. 

We will need an estimate that appeared in \cite[Proposition 2.1]{Bez08}, and  was proved using a Bellman function argument. In that paper and others, a sequence $\{\alpha_I\}_{I\in\mathcal{D}}$ is said to be a $\omega$-Carleson sequence if, and only if, there is a $C>0$ such that for all $J\in\mathcal{D}$ we have that $\sum_{I\in\mathcal{D}(J)} \alpha_I \leq C\,  \omega (J)$. The statement below looks slightly different than what appeared in the original publication because we are  using instead definition~\eqref{w-Carleson-sequence}.


\begin{lem}[Beznosova's Little Lemma \cite{Bez08}]\label{lem:Beznosova} If a sequence of positive numbers $\{\lambda_I\}_{I\in\mathcal{D}}$ is a Carleson sequence, 
then for all weights $\omega $ such that $\omega^{-1}$ is also a weight, we have that the sequence $\{\nu_I\}_{I\in\mathcal{D}}$, defined by $\nu_I=\frac{\lambda_I}{\langle \omega^{-1}\rangle_I \langle \omega\rangle_I }$, is a $\omega$-Carleson sequence.
\end{lem}

\subsection{Sawyer's Estimates}\label{sec:Sawyer}

A key tool needed for us to prove  Theorem~\ref{thm:1}  is  Sawyer's estimate. The following Sawyer's type estimate can be found in \cite[Lemma 9, Proof in Section 6.3]{Per09}, where it is proved using a Bellman's function argument. Although  dyadic doubling is in the hypothesis, doubling of $\sigma$ is not needed in the proof (it is needed in the argument for other Bellman function lemmas on that paper, but not for this one), hence we state the result without that assumption. 

\begin{prop}[$d\sigma$-Sawyer's Estimate]\label{prop:Sawyer's}
Given a positive measure~$\sigma$, a weight $\omega$, a constant $Q>0$, and a sequence of positive numbers $\{\lambda_I \}_{I \in \mathcal{D}}$. Suppose that for all dyadic intervals $J$
\[ \frac{1}{\sigma(J)} \sum_{I \in \mathcal{D}(J)} | \langle \omega \rangle_I^{\sigma}|^2  \lambda_I \leq Q \langle \omega\rangle_J^{\sigma}, \]
then
$$ \frac{1}{\sigma(J)} \sum_{I \in \mathcal{D}(J)} | \langle f \omega^{\frac12}\rangle_I^{\sigma}|^2  \lambda_I \leq 4Q \langle f^2\rangle_J^{\sigma}. $$
\end{prop}
Where here $\langle \omega\rangle_I^{\sigma}$ denotes the integral average of $\omega$ with respect to the measure $\sigma$, that is, $\frac{1}{\sigma(I)} \int_I \omega\,d\sigma$, and $\sigma(I)=\int_Id\sigma$.
%
%

\section{Three Weight Conditions}\label{sec:3weight}

The operators we aim to analyze are inherently tied to a weight. Furthermore, our focus lies on their behavior within the two-weight setting, where distinct weights are assigned to the domain and target spaces. Naturally, any boundedness condition for these operators in such spaces must involve all three weights. 

In a recent study, \cite{CHMPW24}, a three-weight condition is proved to be necessary (also one of four sufficient conditions) for the boundedness of the weighted Haar multipliers \(T_w\), mapping \(L^2(u)\) into \(L^2(v)\). This condition is expressed as:  
\[
\sup_{I \in \mathcal{D}} \frac{\langle vw^2\rangle_I \langle u^{-1}\rangle_I}{\langle w\rangle_I^2} < \infty.
\]  

This criterion can be viewed as a generalization of the classical \(A_p\) condition, which plays a central role in the study of numerous operators within the weighted framework. In what follows, we introduce this extension and discuss various results that establish connections between such conditions.

For weights $u$ and $v$, and  a positive measure $\sigma$, we say that  $(u,v) \in A^d_p(d\sigma)$  if, and only if, 
\begin{equation}
  [u,v]_{A_p^d(d\sigma)}:= \sup_{I \in \mathcal{D}} \langle v\rangle_I^{\sigma} \left ( \langle u^{{-1}/(p-1)}\rangle_I^{\sigma}\right )^{p-1}<\infty.
  \label{ap_geral}
\end{equation}
Note that  for a pair $(u,v)$ to be in joint $A^d_p(d\sigma)$, we need  both $v, u^{\frac{-1}{p-1}} \in L^1_{{\rm loc}}(d\sigma)$. 
We say $u\in A_p^d(d\sigma)$ when $(u,u)\in A_p^d(d\sigma)$.

For a measure $d\sigma = wdx$ where $w$ is a weight, we say that $(u,v) \in A^d_p(w)$  if, and only if, 
\begin{equation}
[u,v]_{A_p^d(w)}:= \sup_{I \in \mathcal{D}}   \langle v\rangle_I^{w} \left ( \langle u^{{-1}/(p-1)}\rangle_I^{w}\right )^{p-1}<\infty.
 \label{ap_weig}
\end{equation}
In this case, for a pair $(u,v)$ to be in joint $A^d_p(w)$, we need $v,  u^{\frac{-1}{p-1}}\in L^1_{{\rm loc}}(w)$, or equivalently $wv, wu^{\frac{-1}{p-1}} \in L^1_{{\rm loc}}(\mathbb{R})$.
In the case $u=v$, we say that $u \in A^d_p(w)$ if, and only if,  $(u,u) \in A^d_p(w)$.

Another important class of weights is the dyadic $p$-Reverse H\"older class, we say that  $ (u,v) \in RH^d_p(w)$ if, and only if, 
\begin{equation}
[u,v]_{RH_p^d(w)}:= \sup_{I \in \mathcal{D}} \left (\langle v^p\rangle_I^w\right )^{1/p}\left (\langle u\rangle_I^w\right )^{-1} <\infty.
    \label{rh_weig}
\end{equation}
In this case,  for a pair $(u,v)$ to be in joint $RH^d_p(w)$, we need $v^p, u\in L^1_{{\rm loc}}(w)$, or equivalently $wv^p, wu \in L^1_{{\rm loc}}(\mathbb{R})$. In the case $u=v$, we say that $u \in RH^d_p(w)$ if, and only if,  $(u,u) \in RH^d_p(w)$.


The following proposition connects many of the previous definitions when $p=2$.

\begin{prop}\label{prop:equiv-3weight-conditions} Given three weights $u,v,w$, so that $vw^2, u^{-1} \in L^1_{{\rm loc}}(\mathbb{R})$,  then the following are equivalent:
{\rm (i)} $ \sup_{I \in \mathcal{D}} \frac{\langle vw^2\rangle_I \langle u^{-1}\rangle_I}{\langle w\rangle_I^2}< \infty$, {\rm (ii)} $ (uw, vw) \in A^d_2(w)$, {\rm (iii)} $ (v^{-1}w^{-1}, u^{-1}w^{-1}) \in A^d_2(w)$, and {\rm (iv)} $(uw,v^{\frac12}u^{\frac12}w) \in RH^d_2(u^{-1})$. 
\end{prop}

\begin{proof} We shall prove that conditions (ii), (iii) and (iv) are equivalent to condition (i). Given three weights $u,v,w$, so that $vw^2, u^{-1} \in L^1_{{\rm loc}}(\mathbb{R})$, then $(uw, vw) \in  A^d_2(w)$ if, and only if,
$$
 \sup_{I \in \mathcal{D}} \langle vw\rangle^w_I  \langle (uw)^{-1}\rangle^w_I = \sup_{I \in \mathcal{D}}\frac{\langle vw^2\rangle_I \langle u^{-1}\rangle_I}{\langle w\rangle_I^2} <\infty.$$
Which proves that condition (i) is equivalent to condition (ii). 

In order to prove that condition (i) is equivalent to condition (iii), note
 that $  (v^{-1}w^{-1}, u^{-1}w^{-1}) \in A^d_2(w) $ if, and only if,
$$\sup_{I \in \mathcal{D}} \langle u^{-1}w^{-1}\rangle_I^w \langle(v^{-1}w^{-1})^{-1}\rangle_I^w = \sup_{I \in \mathcal{D}} \frac{\langle u^{-1}\rangle_I\langle vw^2\rangle_I }{\langle w\rangle_I^2}<\infty. $$




Finally, note that $(uw,v^{\frac12}u^{\frac12}w) \in RH_2(u^{-1})$ if, and only if,
$$ \sup_{I \in \mathcal{D}} {\left \langle \big( v^{\frac12}u^{\frac12}w \big)^2\right \rangle_I^{u^{-1}}}{ \left (\langle wu\rangle_I^{u^{-1}}\right )^{-2}}  = \sup_{I \in \mathcal{D}} \frac{\langle vw^2\rangle_I \langle u^{-1}\rangle_I}{\langle w\rangle_I^{2}}< \infty.$$
Proving that condition (i) is equivalent to condition  (iv).

 \end{proof}

\begin{remark}\label{remark1}
 If $u=v$, then the following are equivalent: {\rm (a)} $\sup_{I \in \mathcal{D}}{\langle uw^2\rangle_I \langle u^{-1}\rangle_I}{\langle w\rangle_I^{-2}}<\infty$, 
 {\rm (b)} $uw \in A^d_2(w)$, {\rm (c)} $u^{-1}w^{-1} \in A^d_2(w)$, and {\rm (d)}  $uw \in RH^d_2(u^{-1})$.
 
  If, in addition, $u=1$, then  the following are equivalent: {\rm (a)} $\sup_{I \in \mathcal{D}}{\langle w^2\rangle_I }{\langle w\rangle_I^{-2}}<\infty$,
 {\rm (b)} $w \in A^d_2(w)$, {\rm (c)} $w^{-1} \in A^d_2(w)$, and {\rm (d)}~$w \in RH^d_2$. 
\end{remark}

We need to introduce  a weak form of the weighted \( A_2(w) \). Similar weaker conditions were introduced in \cite{CaGS15} (for $w=1$), where the authors addressed the limitation that the classical Rubio de Francia extrapolation theorem does not extend to the endpoint \( p = 1 \). To overcome this, they defined a broader class of weights, slightly larger than \( A_p \).
 The broader class introduced by the authors extends the applicability of this extrapolation result. For every \( 1 \leq p < \infty \), they define the restricted \( A_p \) class, \( A_p^R \), is defined as the set of weights \( u \) for which the following quantity is finite:  
\[
\|u\|_{A_p^R} = \sup_{E \subset Q} \left( \frac{|E|}{|Q|} \right) \left( \frac{u(Q)}{u(E)} \right)^{1/p},
\]
\noindent where the supremum is taken over all cubes \( Q \) and all measurable subsets \( E \subset Q \). An application of H\"older's inequality shows that $\|u\|_{A_p^R}\leq [u]_{A_p}^{1/p}$, hence $A_p\subset A_p^R$.

We need to extend this definition in two directions: we need a two weight type $A_2$ condition and also allowing an underlying measure $w$. Our condition does not need to hold for every descendant of an interval. Instead, we can slightly weaken it by considering only the ratio between an interval and its immediate children. We say that \((u, v)\) belongs to the \emph{restricted  $\tilde{A}^d_2(w)$ class} if  
\begin{equation}\label{def:restrictedA2}
[u,v]_{\tilde{A}^d_2(w)}:=\sup_{J \in \mathcal{D}} \left(\frac{w(J^*)}{w(\widetilde{J})}\right)^2 \langle u^{-1} \rangle^w_J \langle v \rangle^w_J <\infty,
\end{equation} 
where  as before, \( J^* \) and \( \widetilde{J} \) are the sibling and the parent of \(J\) respectively.

Note that for any weight $w$, we have that  \({w(J^*)}/{w(\widetilde{J})} \leq 1 \), thus for any pair of weights~\( (u,v) \), we  see that,   \([u,v]_{\tilde{A}^d_2(w)} \leq [u,v]_{A^d_2(w)} \), which implies that \(A^d_2(w) \subset \tilde{A}^d_2(w) \). Moreover, if \( w \) is dyadic doubling then \( \tilde{A}^d_2(w) = A_2(w) \).


\section{Weighted Square Function}\label{SquareFunction}\label{sec:Sw}

In this section, we will prove the main results about the weighted square function, aiming to establish its boundedness in the two-weighted setting. Specifically, we will prove the necessity and sufficiency of the conditions introduced in Theorems~\ref{thm:1} and~\ref{thm:2}. First, we will show some sufficient conditions for the boundedness of the weighted dyadic square function that do not require any additional assumptions on the weights. We can prove that these conditions are necessary under the additional hypothesis that $w$ is a dyadic doubling weight (Corollary~\ref{coro:iff}). In Theorem~\ref{thm:2}  (restated as Theorem~\ref{thm:necessity} in this section) we find necessary conditions for the boundedness of the weighted dyadic square function, which almost match the sufficient conditions in Theorem~\ref{thm:1} (restated as  Theorem~\ref{thm:square-function-two-weight} in this section), except that we need to weaken the three-weight condition.

\begin{thm}[The sufficiency]\label{thm:square-function-two-weight} Let $u,v,w$ be weights, such that  $u^{-1}w^2\in L^1_{{\rm loc}}(\mathbb{R})$, and 
\begin{itemize}
    \item[{\rm (i)}] {\rm (Three-weight condition)} $(uw^{-1},vw^{-1}) \in A^d_2(w)$,
     \item[{\rm (ii)}] {\rm (Carleson condition)}  there is a $C>0$ such that for all $J\in\mathcal{D}$ we have that,
$$\frac{1}{w(J)} \sum_{I \in \mathcal{D}(J) } \big| \Delta^w_I (u^{-1}w) \big|^2  \left [ \frac{w(I^-)^2v(I^+) + w(I^+)^2v(I^-)}{w(I)^2} \right ]  \leq C \langle u^{-1}w\rangle_J^w.$$
\end{itemize}
Then, the weighted dyadic square function $S^d_w$
is bounded from $L^2(u)$ into $L^2(v)$. 
Furthermore,
\[\|S_w^d\|_{L^2(u)\to L^2(v)} \leq \sqrt{[uw^{-1},vw^{-1}]_{A_2^d(w)}}+ 2\sqrt{C}.\]
\end{thm}


\begin{proof}[Proof of Theorem~\ref{thm:square-function-two-weight}]
First observe that we can write $\|S_w^d f \|^2_{L^2(v)}$ as
$$\|S_w^d f \|^2_{L^2(v)} =\sum_{I \in \mathcal{D}} |\langle f\rangle^w_I  - \langle f\rangle^w_{\widetilde{I}}|^2 v(I).$$
Therefore by~\eqref{eq:Delta^w}, replacing $ |\langle f \rangle^w_I  - \langle f\rangle^w_{\widetilde{I}} |^2 =\frac{w(I^*)}{w(\widetilde{I})w(I)}\big| \langle f, h^w_{\widetilde{I}} \rangle_w \big|^2 $, and  changing summation over children to summation over parents, we get
\begin{align*}
 \|S_w^d f \|^2_{L^2(v)} &=\sum_{I \in \mathcal{D}} \frac{1}{w(I) }\Bigg( \frac{v(I^+)w(I^-)}{w(I^+)} + \frac{v(I^-)w(I^+)}{w(I^-)} \Bigg)  \big| \langle f, h^w_{I} \rangle_w \big|^2   \\
 &=\sum_{I \in \mathcal{D}} \Bigg( \frac{w(I^-)\langle vw^{-1}\rangle_{I^+}^w+w(I^+)\langle vw^{-1}\rangle_{I^-}^w}{w(I)} \Bigg)  \big| \langle f, h^w_{I} \rangle_w \big|^2\\
   &= \sum_{I \in \mathcal{D}} K_I^{w,v} \big| \langle f, h^w_{{I}} \rangle_w \big|^2,
\end{align*}
where on the second equality, all we have done is write $v(J)/w(J)=\langle vw^{-1}\rangle_J^w$, and where 
\[ K_I^{w,v}:= \frac{w(I^-)\langle vw^{-1}\rangle_{I^+}^w +w(I^+)\langle vw^{-1}\rangle_{I^-}^w}{w(I) }. \]
Boundedness of $K_I^{w,v}$ would be the necessary and sufficient condition for boundedness of $S_w$ from $L^2(w)$ into $L^2(v)$. For the necessity, by testing on $f=h_I^w$, and for the sufficiency,  by Bessel's inequality applied to the orthonormal system $\{h_I^w\}_{I\in\mathcal{D}}$ in $L^2(w)$. However, we are interested in boundedness from $L^2(u)$ into $L^2(v)$, so we will replace the Haar functions with respect to $w$, by Haar functions with respect to $w^2u^{-1}$, using~\eqref{eqn:weightedHaar-linear-combination}, namely,
\begin{equation}\label{haarsystem}
 h_I^w  = \alpha_I^{u,w}h_I^{w^2u^{-1}} +\beta_I^{u,w}\mathbbm{1}_I, 
 \end{equation}
 where,
 \begin{eqnarray}
 \alpha_I^{u,w}:= \sqrt{\frac{\langle wu^{-1}\rangle^w_{I^+} \langle wu^{-1}\rangle^w_{I^-}}{\langle wu^{-1}\rangle^w_{I}}}, \label{alpha}\\
 \beta_I^{u,w}=\sqrt{\frac{w(I^+) w(I^-)}{w(I)} } \frac{1}{w(I)} \frac{\Delta^w_I (wu^{-1})}{\langle wu^{-1}\rangle_I^w}. \label{beta}
\end{eqnarray}

Also remember that $f \in L^2(w)$ if and only if $u^{-\frac12}w^{\frac12} f \in L^2(u)$, moreover 
$\|u^{-\frac12}w^{\frac12} f\|_{L^2(u)}=\|f\|_{L^2(w)}$. Then, we want to bound for all $ f \in L^2(w)$, 
\[ \left \|S_w^d (u^{-\frac12}w^{\frac12} f) \right \|^2_{L^2(v)} = \sum_{I \in \mathcal{D}} K_I^{w,v} \left | \langle u^{-\frac12}w^{\frac12} f, h^w_{{I}} \rangle_w \right |^2 \leq C\|u^{-\frac12}w^{\frac12} f \|^2_{L^2(u)}= C \| f\|^2_{L^2(w)}.  \]
Substituting \eqref{haarsystem} in $ \displaystyle \sum_{I \in \mathcal{D}} K_I^{w,v} \big| \langle u^{-\frac12}w^{\frac12} f, h^w_{{I}} \rangle_w \big|^2 $, we have that
\begin{align*}
 \left \|S_w^d (u^{-\frac12}w^{\frac12} f) \right \|^2_{L^2(v)} &=  \sum_{I \in \mathcal{D}} K_I^{w,v} \left | \langle u^{-\frac12}w^{\frac12} f,\alpha_I^{u,w}h_I^{w^2u^{-1}} +\beta_I^{u,w}\mathbbm{1}_I \rangle_w \right |^2 \\
 &=  \sum_{I \in \mathcal{D}} K_I^{w,v} \left | \alpha_I^{u,w}\langle u^{\frac12}w^{-\frac12} f,h_I^{w^2u^{-1}}\rangle_{w^2u^{-1} } +
 \beta_I^{u,w} \langle u^{\frac12}w^{-\frac12} f , \mathbbm{1}_I \rangle_{w^2u^{-1} }\right |^2 \\
 &= \Sigma_1 + \Sigma_2 + \Sigma_3.
\end{align*}
In the second equality, we used that $u^{-1/2}w^{1/2}w=u^{1/2}w^{-1/2}w^2u^{-1}$, to move some weights to the inner product.
In the third equality, after squaring we get three terms ($(a+b)^2=a^2+2ab+b^2$), defined by:
\begin{align}
\Sigma_1: &= \sum_{I \in \mathcal{D}} K_I^{w,v} \big|\alpha_I^{u,w} \big|^2 \left  | \langle u^{\frac12}w^{-\frac12} f,h_I^{w^2u^{-1}}\rangle_{w^2u^{-1} }\right |^2, \label{sigma_1}\\
\Sigma_2: &= \sum_{I \in \mathcal{D}} K_I^{w,v} 2 \alpha_I^{u,w} \beta_I^{u,w}  \langle u^{\frac12}w^{-\frac12} f,h_I^{w^2u^{-1}}\rangle_{w^2u^{-1} } \langle u^{\frac12}w^{-\frac12} f , \mathbbm{1}_I \rangle_{w^2u^{-1} },\label{sigma_2} \\
\Sigma_3:&=\sum_{I \in \mathcal{D}} K_I^{w,v} \big | \beta_I^{u,w} \big|^2 \left | \langle u^{\frac12}w^{-\frac12} f , \mathbbm{1}_I \rangle_{w^2u^{-1} }\right |^2. \label{sigma_3}
\end{align}
It is enough to bound $\Sigma_1$ and $\Sigma_3$, since, by the Cauchy-Schwarz inequality, $|\Sigma_2 | \leq 2\sqrt{\Sigma_1 \Sigma_3}$.\\

\noindent{\bf Bound for $\Sigma_1$:} Substituting \eqref{alpha} in \eqref{sigma_1}, we have
\begin{align}
\Sigma_1: &= \sum_{I \in \mathcal{D}} K_I^{w,v} \frac{\langle wu^{-1}\rangle^w_{I^+} \langle wu^{-1}\rangle^w_{I^-}}{\langle wu^{-1}\rangle^w_{I}}
 \left | \langle u^{\frac12}w^{-\frac12} f,h_I^{w^2u^{-1}}\rangle_{w^2u^{-1} }\right |^2, \label{Sigma1_estimate} 
\end{align}
where
\begin{align*}
 K_I^{w,v} \frac{\langle wu^{-1}\rangle^w_{I^+} \langle wu^{-1}\rangle^w_{I^-}}{\langle wu^{-1}\rangle^w_{I}}&= \frac{\big (w(I^-)\langle vw^{-1}\rangle^w_{I^+} +w(I^+)\langle vw^{-1}\rangle^w_{I^-} \big )}{w(I)}  \frac{\langle wu^{-1}\rangle^w_{I^+}  \langle wu^{-1}\rangle^w_{I^-} }{ \langle wu^{-1}\rangle^w_{I}  }\nonumber \\
 &\hskip -.5in=\frac{w(I^-)\langle vw^{-1}\rangle^w_{I^+} \langle wu^{-1}\rangle^w_{I^+} \langle wu^{-1}\rangle^w_{I^-} + w(I^+)\langle vw^{-1}\rangle^w_{I^-} \langle wu^{-1}\rangle^w_{I^+} \langle wu^{-1}\rangle^w_{I^-} }{w(I)  \langle wu^{-1}\rangle^w_{I}  }  \nonumber \\
 & \hskip -.5in\leq [uw^{-1},vw^{-1}]_{A^d_2(w)}\frac{w(I^-)\langle wu^{-1}\rangle^w_{I^-} + w(I^+)\langle wu^{-1}\rangle^w_{I^+} }{w(I)  \langle wu^{-1}\rangle^w_{I}  }.  \nonumber
\end{align*}
The last inequality follows because $ [uw^{-1},vw^{-1}]_{A^d_2(w)} \geq \langle vw^{-1}\rangle^w_{J} \langle wu^{-1}\rangle^w_{J}$ for all  $J \in \mathcal{D}$. Thus, observing that the numerator equals the denominator in the last estimate, that is, $w(I^-)\langle wu^{-1}\rangle^w_{I^-} + w(I^+)\langle wu^{-1}\rangle^w_{I^+} =w(I)\langle wu^{-1}\rangle^w_{I}$, we conclude that,
\begin{align}
 K_I^{w,v} \frac{\langle wu^{-1}\rangle^w_{I^+} \langle wu^{-1}\rangle^w_{I^-}}{\langle wu^{-1}\rangle^w_{I}} \leq 
 [uw^{-1},vw^{-1}]_{A^d_2(w)}. \label{N:estimate}
\end{align}
Therefore, substituting \eqref{N:estimate} in \eqref{Sigma1_estimate} and using Bessel's inequality we have:
\begin{align*}
\Sigma_1 &\leq [uw^{-1},vw^{-1}]_{A^d_2(w)} \sum_{I \in \mathcal{D}}  \Big| \big\langle u^{\frac12}w^{-\frac12} f,h_I^{w^2u^{-1}}\big\rangle_{w^2u^{-1} }\Big|^2 \\
 &\leq [uw^{-1},vw^{-1}]_{A^d_2(w)} \big\| u^{\frac12}w^{-\frac12} f \big\|^2_{L^2(w^2u^{-1})} \, = \,[uw^{-1},vw^{-1}]_{A^d_2(w)}\big\| f \big\|^2_{L^2(w)}.
\end{align*}

\noindent{\bf Bound for $\Sigma_3$:} Substituting \eqref{beta} in \eqref{sigma_3}, we have
\begin{align*}
\Sigma_3 &= \sum_{I \in \mathcal{D}} K_I^{w,v} \big| \beta_I^{u,w} \big|^2 \left | \langle u^{\frac12}w^{-\frac12} f , \mathbbm{1}_I \rangle_{w^2u^{-1} }\right|^2 \\
&= \sum_{I \in \mathcal{D}} K_I^{w,v} \frac{w(I^+)w(I^-)}{w(I)} \Bigg[ \frac{\Delta_I^w(wu^{-1}) }{\langle wu^{-1}\rangle^w_{I} }\Bigg]^2 \left | \frac{\langle u^{\frac12}w^{-\frac12} f , \mathbbm{1}_I \rangle_{w^2u^{-1} }}{w(I)}\right |^2 \\
&= \sum_{I \in \mathcal{D}} K_I^{w,v} \frac{w(I^+)w(I^-)}{w(I)}  \frac{\big|\Delta^w_{I}(wu^{-1})\big|^2}{\big( \langle wu^{-1}\rangle^w_{I} \big)^2} \Big(\langle u^{-\frac12}w^{\frac12} f \rangle_I^w\Big)^2.
\end{align*}
Besides, since $w(J)\langle vw^{-1}\rangle_J^w=v(J)$ for all $J\in\mathcal{D}$, 
\begin{align*}
    K_I^{w,v} \frac{w(I^+)w(I^-)}{w(I)}&= \frac{\big [w(I^-)\langle vw^{-1}\rangle^w_{I^+} +w(I^+) \langle vw^{-1}\rangle^w_{I^-} \big ]}{w(I) } \, \frac{w(I^+)w(I^-)}{w(I)} \\
    &= \frac{w(I^-)^2w(I^+) \langle vw^{-1}\rangle^w_{I^+} +w(I^+)^2w(I^-) \langle vw^{-1}\rangle^w_{I^-} }{w(I)^{2}}  \\
&=\frac{w(I^-)^2v(I^+)+w(I^+)^2v(I^-)}{w(I)^{2}}. 
\end{align*}
Then 
\[ \Sigma_3 = \sum_{I \in \mathcal{D}}  \left [ \frac{w(I^-)^2v(I^+) + w(I^+)^2v(I^-)}{w(I)^2} \right ]  \frac{\big|\Delta^w_{I}(wu^{-1})\big|^2}{\big( \langle wu^{-1}\rangle^w_{I} \big)^2} \Big(\langle u^{-\frac12}w^{\frac12} f \rangle_I^w\Big)^2.\]
Let 
$$\displaystyle{ \lambda_I=  \left [ \frac{w(I^-)^2v(I^+) + w(I^+)^2v(I^-)}{w(I)^2} \right ]  \frac{\big|\Delta^w_{I}(wu^{-1})\big|^2}{\big( \langle wu^{-1}\rangle^w_{I} \big)^2}}, $$  
then by hypothesis (ii),
\[ \frac{1}{w(J)} \sum_{I \in \mathcal{D}(J)} \big( \langle wu^{-1}\rangle^w_{I} \big)^2 \lambda_I \leq C \langle wu^{-1}\rangle_J^w,
 \]
 which allows us to use Sawyer's Estimate to complete the proof. 
  Indeed, Sawyer's Estimate (Proposition~\ref{prop:Sawyer's}) gives us that for all $J\in\mathcal{D}$,
 \[ \frac{1}{w(J)}\sum_{I\in \mathcal{D}(J)} \big (\langle u^{-\frac12}w^{\frac12} f\rangle_I^w \big )^2 \lambda_I \leq 4C \langle f^2\rangle_J^w,\]
 which  implies that for all $J\in\mathcal{D}$, 
 \[ \sum_{I\in \mathcal{D}(J)} \big (\langle u^{-\frac12}w^{\frac12} f\rangle_I^w \big )^2 \lambda_I \leq 4C \int_Jf^2(x)\,w(x)\, dx. \] 
 Adding these estimates for $J^+_n=[0,2^n)$ and for $J^-_{n}=[-2^n,0)$, and letting $n\to\infty$, implies that
 \[ \Sigma_3\leq  \sum_{I\in\mathcal{D}}\big (\langle u^{-\frac12}w^{\frac12} f\rangle_I^w \big )^2 \lambda_I  \leq 4C \|f\|^2_{L^2(w)}. \]
All together, these estimates imply that, 
\[ \|S_w^d\|^2_{L^2(u)\to L^2(v)} \leq 
\big ( \sqrt{[uw^{-1},vw^{-1}]_{A_2(w)^d }} + 2\sqrt{C}\big )^2.\]


\end{proof}


We cannot prove that condition (i) is necessary without assuming that the weight \( w \) is doubling. However, we can establish that a weak form of the weighted \( A_2(w) \) condition is necessary. Namely, the restricted $\tilde{A}_2^d(w)$ condition introduced in Section~\ref{sec:3weight}, is a necessary condition for the two weight boundedness of the weighted dyadic square function. Recall that  \(A^d_2(w) \subset \tilde{A}^d_2(w) \). Moreover, if \( w \) is dyadic doubling then \( \tilde{A}^d_2(w) = A_2(w) \).

\begin{thm}[The necessity]\label{thm:necessity}
 Let $u,v,w$ be weights, such that  $u^{-1}w^2,v\in L^1_{{\rm loc}}(\mathbb{R})$. Suppose $S_w$ is bounded from $L^2(u)$ into $L^2(v)$, then   
\begin{itemize}
    \item[{\rm (i')}] $(uw^{-1},vw^{-1}) \in \tilde{A}^d_2(w)$,
     \item[{\rm (ii)}]  there is a constant $C>0$ such that for all $J\in\mathcal{D}$ we have that
$$\frac{1}{w(J)} \sum_{I \in \mathcal{D}(J) } \big| \Delta^w_I (u^{-1}w) \big|^2  \left [ \frac{w(I^-)^2v(I^+) + w(I^+)^2v(I^-)}{w(I)^2} \right ]  \leq C \langle u^{-1}w\rangle_J^w.$$
\end{itemize}
Moreover $ [uw^{-1},vw^{-1}]_{\tilde{A}_2^d(w)}\leq \|S_w^d\|^2_{L^2(u)\to L^2(v)}$ and $C\leq \|S_w^d\|^2_{L^2(u)\to L^2(v)}$.
\end{thm}
\begin{proof}To get the three-weight condition (i') we will essentially mimic the argument in \cite[Theorem in page 921]{NTV99}. We will test on the functions $f=u^{-1}w\mathbbm{1}_J$, note that 
$\|u^{-1}w\mathbbm{1}_J\|_{L^2(u)}^2= \int_J u^{-1}(x) w^2(x)\, dx<\infty $ by assumption, so $f\in L^2(u)$. From the definition of the weighted square function,
\[ S^d_wf(x)= \bigg (\sum_{I\in\mathcal{D}} |\langle f\rangle_I^w-\langle f\rangle_{\widetilde{I}}^w|^2 \mathbbm{1}_I(x)\bigg )^{1/2},\]
we can estimate from below $S_w(u^{-1}w\mathbbm{1}_J)(x)$ for each $x\in J$. In fact, since $\int_{\widetilde{J}} h\mathbbm{1}_J=\int_Jh$, 
\begin{eqnarray*}
S^d_w(u^{-1}w\mathbbm{1}_J)(x) & \geq & \left | \langle u^{-1}w\mathbbm{1}_J\rangle_J^w -\langle u^{-1}w\mathbbm{1}_J\rangle^w_{\widetilde{J}}\right |\\
& = & \left | \frac{\int_J u^{-1}(x)\,w^2(x)\, dx}{\int_J w(x)\, dx}-\frac{\int_J u^{-1}(x)\,w^2(x)\, dx}{\int_{\widetilde{J}} w(x)\, dx}\right |\\
&= & \frac{\int_J u^{-1}(x)\,w^2(x)\, dx}{\int_J w(x)\, dx} \frac{w(J^*)}{w(\widetilde{J})},
\end{eqnarray*}
where $J^*$ is the sibling of $J$.
We can now estimate the $L^2(v)$ norm of $S^d_w(f)$ from below, and use the hypothesis to conclude that
\[\frac{\big (\int_J u^{-1}(x)\,w^2(x)\, dx\big )^2}{\big (\int_J w(x)\, dx\big )^2} \left (\frac{w(J^*)}{w(\widetilde{J})}\right )^2 v(J) \leq \|S_w^d(f)\|_{L^2(v)}^2 \leq C \int_{J} u^{-1}(x)\, w^2(x)\,dx.\]
Canceling, rearranging, and multiplying and dividing by $|J|^2$, we conclude that 
\[\frac{\langle u^{-1}w^2\rangle_J \langle v\rangle_J}{\langle w\rangle_J^2} \left (\frac{w(J^*)}{w(\widetilde{J})}\right )^2= \frac{\int_J u^{-1}(x)\,w^2(x)\, dx}{(w(J))^2} \left (\frac{w(J^*)}{w(\widetilde{J})}\right )^2 v(J) \leq C.\]
Hence $(u^{-1}w,vw^{-1})\in \tilde{A}_2^d(w)$ and  condition (i') holds, with $C= \|S_w^d\|^2_{L^2(u)\to L^2(v)}$.

As before, let $f=u^{-1}w\mathbbm{1}_J$, since $u^{-1}w^2 \in L^1_{{\rm loc}}(\R)$ then  $f\in L^2(u)$ and
\[\|f\|_{L^2(u)}^2= \int_J u^{-2} w^2 u\, dx = w(J)\, \langle  u^{-1}w\rangle_J^w.\]
Moreover,
\begin{eqnarray}
\|S_w(f)\|_{L^2(v)}^2 & = & \sum_{I \in \mathcal{D}} \left |\langle  u^{-1}w\mathbbm{1}_J\rangle_I^w-  \langle u^{-1}w\mathbbm{1}_J\rangle^w_{\widetilde{I}}\right |^2  v(I) \nonumber\\
&\geq &  \sum_{I \in \mathcal{D}:I\subsetneq J}   \left | \langle u^{-1}w\rangle_I^w- \langle u^{-1}w\rangle^w_{\widetilde{I}}\right |^2  v(I). \quad\quad   \label{eq:lowerboundSquare}
\end{eqnarray}

We need to show that condition (ii) holds, that is,
\begin{equation}\label{def:Sigma3}
\Sigma_3:= \sum_{I \in \mathcal{D}(J)}   \left | \Delta_I^w(u^{-1}w)\right |^2 \left [ \frac{w(I^-)^2v(I^+) + w(I^+)^2v(I^-)}{w(I)^2} \right ] \leq C  w(J) \langle  u^{-1}w\rangle_J^w.
\end{equation}
It will suffice to show that
\begin{equation}\label{eq:comparison-weighted-differences2}
\Sigma_3 = \sum_{I \in \mathcal{D}:I\subsetneq J}   \left | \langle u^{-1}w\rangle_I^w- \langle u^{-1}w\rangle^w_{\widetilde{I}}\right |^2 v(I). 
\end{equation}
If~\eqref{eq:comparison-weighted-differences2} holds,  then, since we are assuming that $S^d_w$ is bounded from $L^2(u)$ into $L^2(v)$, we will conclude that,
\[ \Sigma_3\leq \|S_w^d(u^{-1}w\mathbbm{1}_J)\|_{L^2(v)}^2  \leq C w(J)  \langle  u^{-1}w\rangle_J^w,\]
which is condition (ii), with $C=\|S_w^d\|_{L^2(u)\to L^2(v)}$.

To verify \eqref{eq:comparison-weighted-differences2}, let us use \eqref{eq:DeltaIw(f)}, namely, substitute $|\Delta_I^w(f)|^2=\frac{w(I)}{w(I^-)w(I^+)} |\langle f,h_I^w\rangle_w|^2$ into \eqref{def:Sigma3}, the definition of $\Sigma_3$,  and simplify,  to get,
\begin{eqnarray*}
\Sigma_3  & = &
\sum_{I \in \mathcal{D}(J)}  \left | \langle u^{-1}w,h_I^w\rangle_w\right |^2 \left [ \frac{w(I^-)v(I^+)}{w(I)w(I^+)} + \frac{w(I^+)v(I^-)}{w(I)w(I^-)} \right ].
\end{eqnarray*}
Next, let us first use \eqref{eq:Delta^w} for $I^{\pm}$, namely, $|\langle f,h_I^w\rangle_w|^2=\frac{w(I)w(I^{\pm})}{w(I^{\mp})}|\langle f\rangle_{I^{\pm}}^w-\langle f\rangle_I^w|^2$,  and simplify, then switch the summation index from parent to children, to get,
\begin{eqnarray*} 
\Sigma_3 & = &  \sum_{I \in \mathcal{D}(J)} \left ( \left | \langle u^{-1}w\rangle_{I^+}^w -\langle u^{-1}w\rangle_{I}^w\right |^2v(I^+) + \left | \langle u^{-1}w\rangle_{I^-}^w -\langle u^{-1}w\rangle_{I}^w\right |^2 v(I^-)\right ) \\
& = &  \sum_{I \in \mathcal{D}:I\subsetneq J}  \left | \langle u^{-1}w\rangle_I^w- \langle u^{-1}w\rangle^w_{\widetilde{I}}\right |^2 v(I).
\end{eqnarray*}

\end{proof}

If the weight $w$ is doubling, then condition (i) in Theorem~\ref{thm:square-function-two-weight} and condition (i') in Theorem~\ref{thm:necessity} are the same, which allows us to write the following corollary:

\begin{coro}\label{coro:iff}
 Let $u,v,w$ be weights, such that $w$ is doubling,  $u^{-1}w^2,v\in L^1_{{\rm loc}}(\mathbb{R})$. Then, $S_w$ is bounded from $L^2(u)$ into $L^2(v)$ if, and only if,
\begin{itemize}
    \item[{\rm (i)}] $(uw^{-1}, vw^{-1}) \in A^d_2(w)$,
     \item[{\rm (ii)}]  there is a constant $C>0$ such that for all $J\in\mathcal{D}$ we have that
$$\frac{1}{w(J)} \sum_{I \in \mathcal{D}(J) } \big| \Delta^w_I (u^{-1}w) \big|^2  \left [ \frac{w(I^-)^2v(I^+) + w(I^+)^2v(I^-)}{w(I)^2} \right ]  \leq C \langle u^{-1}w\rangle_J^w.$$
\end{itemize}
\end{coro}

\section{Previous Known Results are Recovered}\label{sec:previous}

In this section, we will show that the results of Section~\ref{sec:Sw} recover many known results about the square function. 
%
%
%

For our first result we will need Beznosova's Little Lemma (Lemma~\ref{lem:Beznosova} introduced in Section~\ref{sec:Carleson}.) Recall that   if $u\in A_2^d$ or, equivalently  $u^{-1}\in A_2^d$, then $u^{-1}\in A_{\infty}^d$, and the sequence $\lambda_I=|\Delta_I (u^{-1})|^2 / \langle u^{-1}\rangle_I^2$ is a Carleson sequence with intensity $B\lesssim \log ([u^{-1}]_{A_{\infty}^d})$, see Proposition~\ref{thm:FKP} and \cite{FKP91,BR14}. We can apply Lemma~\ref{lem:Beznosova}, seting $\omega=u^{-1}$,  to get a  $u^{-1}$-Carleson sequence ``for free", namely the sequence $\frac{\lambda_I}{\langle u\rangle_I \langle u^{-1}\rangle_I}$. In other words it is true that for all $J\in\mathcal{D}$,
 \begin{equation}\label{Little-Lemma}
 \frac{1}{|J|} \sum_{I\in \mathcal{D}(J)} |I| \left |\frac{\Delta_I u^{-1}}{\langle u^{-1}\rangle_I}\right |^2 \frac{1}{ \langle u\rangle_I} \lesssim \log([u^{-1}]_{A_{\infty}^d}) \langle u^{-1}\rangle_J.
 \end{equation}

To state our first corollary, we set $w=1$ and $u=v$ in Theorem~\ref{thm:1}, to obtain the well known one-weight result for the unweighted dyadic square function, S. Buckley presented a different elementary argument \cite[Section 3, Theorem 3.6]{Buc93}. There are one-weight results for the Walsh-Paley system on $[0,1]$,  that go back to R. Gundy and R. Wheeden \cite[Theorem 2]{GW74}, where they prove the qualitative result using good-lambda inequalities. 

\begin{coro}[\cite{GW74, Buc93}]\label{cor:Wilson}
Given a weight $u$, then the dyadic square function $S^d$ is bounded on $L^2(u)$ if, and only if, $ u \in A_2^d$.
\end{coro}
\begin{proof}
In this case, we have $v=u$ and $w=1$, which implies that condition (i)  in Theorem~\ref{thm:1}  will become:
$ (u,u) \in A_{2}^d$ which is equivalent to  $u \in A_2^d$.

On the other hand,  condition (ii) reads:
$\displaystyle  \frac{1}{|J|} \sum_{I \in \mathcal{D}(J) } |I| \big| \Delta_I (u^{-1}) \big|^2 \langle u\rangle_I   \leq C \langle u^{-1}\rangle_J$.

Since $u\in A_2^d$, we  can now use~\eqref{Little-Lemma}  to get condition (ii), namely,
\begin{eqnarray*}
 \frac{1}{|J|} \sum_{I \in \mathcal{D}(J) } |I| \big| \Delta_I (u^{-1}) \big|^2 \langle u\rangle_I  & =&
\frac{1}{|J|} \sum_{I \in \mathcal{D}(J) } |I| \left |\frac{ \Delta_I (u^{-1})}{\langle u^{-1}\rangle_I}\right |^2  \frac{\big (\langle u^{-1}\rangle_I \langle u\rangle_I \big )^{2}}{\langle u\rangle_I}\\
& \leq &
[u]_{A_2^d}^2 \frac{1}{|J|} \sum_{I\in \mathcal{D}(J)} |I| \left |\frac{\Delta_I u^{-1}}{\langle u^{-1}\rangle_I}\right |^2 \frac{1}{ \langle u\rangle_I} \\
& \lesssim &  \log ([u^{-1}]_{A_{\infty}^d})\, [u]_{A_2^d}^2  \langle u^{-1}\rangle_J.
\end{eqnarray*}
 \end{proof}
This argument does not recover the optimal linear bound on the $A_2$-characteristic~\cite{HTV00, Wit00}. Theorem~\ref{thm:1}, and the proof of Corollary~\ref{cor:Wilson}, provide  a bound that is off the linear bound by the square root of a logarithmic factor, more precisely, 
\[ \|S^d\|_{L^2(u)\to L^2(u)} \lesssim [u]_{A_2^d} \big (\log([u^{-1}]_{A_{\infty}^d})\big )^{1/2}.\]

The second corollary pertains the two-weight boundedness of the unweighted dyadic square function.
In \cite[Theorem in p. 921]{NTV99}, a Sawyer-type theorem is proved that gives testing conditions as necessary and sufficient for the two-weight boundedness of the (generalized) dyadic square function. They reduce the sufficiency to verifying the conditions listed in Corollary~\ref{coro:NTV}. They do not state that the conditions in Corollary~\ref{coro:NTV} are necessary but this is obvious from their argument. We record their theorem for the dyadic square function.

\begin{thm}[\cite{NTV99}]
Given weights $u,v$, the dyadic square function $S^d$ is bounded from $L^2(u)$ into $L^2(v)$ if and only if  there is a constant $C>0$ such that 
\begin{equation}\label{testing}
 \int_I |S^d(u^{-1}\mathbbm{1}_I)(x)|^2 v(x)\,dx \leq C \int_I u^{-1}(x) \,dx, \quad\quad\mbox{for all $I\in\mathcal{D}$}.
 \end{equation}
\end{thm}
Note that \eqref{testing} are localized testing conditions on the functions $f=u^{-1}\mathbbm{1}_I$.

\begin{coro}[\cite{NTV99}]\label{coro:NTV}
Given weight $u,v$, then the dyadic square function $S^d$ bounded from $L^2(u)$ into $L^2(v)$ if and only if
    \begin{itemize}
        \item[{\rm (i)}] $(u,v) \in A^d_2$.
        \item[{\rm (ii)}] There is a constant $C>0$ such that for all $J\in \mathcal{D}$ we have 
        \[\frac{1}{|J|} \sum_{I\in\mathcal{D}(J)}|I|| \Delta_I u^{-1} |^2 \langle v\rangle_I \leq C\langle u^{-1}\rangle_J.\]
    \end{itemize}
\end{coro}
\begin{proof}
    Consider $w=1$ in conditions (i) and (ii) of Corollary~\ref{coro:iff} and we will have exactly the  conditions in this corollary.
\end{proof}

The third  corollary pertains the unweighted bound for the weighted square function. We will need a couple auxiliary results. The first is Buckley's summation characterization of $RH_2^d$ (Proposition~\ref{prop:Buckley}). The second is Proposition~\ref{lem:CHMPW} \cite[Proposition 1]{CHMPW24}.
%

\begin{coro}\label{coro:Per}
Given a weight $w$, then the weighted dyadic square function $S^d_w$ is bounded on $L^2(\mathbb{R})$ if, and only if, $ w \in RH^d_2$.
\end{coro}
\begin{proof}
  Now, we have $v=u=1$, then condition (i) in Theorem~\ref{thm:1} will become:
$ (w^{-1},w^{-1}) \in A_{2}^d(w)$ if, and only if, $w^{-1} \in A_2^d(w)$,
  which is equivalent to $w\in RH_2^d$ by Proposition~\ref{prop:equiv-3weight-conditions} (see Remark~\ref{remark1}).
  
In this case, condition (ii) reads 
\begin{equation}\label{eq:RH2}
\Sigma_1:=\frac{1}{|J|}\sum_{I\in\mathcal{D}(J)} |I| \left | \frac{\langle w^2\rangle_{I_+}}{\langle w\rangle_{I_+}} - \frac{\langle w^2\rangle_{I_-}}{\langle w\rangle_{I_-}}\right |^2\leq C \langle w^2\rangle_J.
\end{equation}
We will show that if $w\in RH_2^d$ then \eqref{eq:RH2} holds. Indeed, starting with the left-hand-side of \eqref{eq:RH2}, adding and subtracting $\langle w^2\rangle_{I_{\pm}} / \langle w\rangle_I$ and using that $|a+b+c|^2\leq 2(a^2+b^2+c^2)$, we can estimate $\Sigma_1$ with three sums. Namely,
\begin{eqnarray*}
\Sigma_1& \lesssim &\frac{1}{|J|}\sum_{I\in\mathcal{D}(J)} |I| \langle w^2\rangle_{I_+}^2  \left ( \frac{1}{\langle w\rangle_{I_+}}-\frac{1}{\langle w\rangle_I}\right )^2 
 +\frac{1}{|J|}\sum_{I\in\mathcal{D}(J)} |I| \left | \frac{\langle w^2\rangle_{I_+}-\langle w^2\rangle_{I_-}}{\langle w\rangle_{I}}\right |^2\\
& & \quad\quad\quad +\, \frac{1}{|J|}\sum_{I\in\mathcal{D}(J)} |I| \langle w^2\rangle_{I_-}^2  \left ( \frac{1}{\langle w\rangle_{I_-}}-\frac{1}{\langle w\rangle_I}\right )^2. 
\end{eqnarray*}
We can put together the first and  last  sums,  use that $|\Delta_I(w)|=\frac12|\langle w\rangle_{I_{\pm}}-\langle w\rangle_I|$, then use the $RH_2^d$ property to get that,
\begin{eqnarray*}
\Sigma_1 & \lesssim & 
\frac{1}{|J|}\sum_{I\in\mathcal{D}(J)} |I| \left (\frac{\langle w^2\rangle_{I_+}^2}{\langle w\rangle_{I_+}^2} + \frac{\langle w^2\rangle_{I_-}^2}{\langle w\rangle_{I_-}^2}\right )\left |\frac{\Delta_Iw}{\langle w\rangle_I}\right |^2 +
\frac{1}{|J|}\sum_{I\in\mathcal{D}(J)} |I| \left | \frac{\Delta_I(w^2)}{\langle w\rangle_{I}}\right |^2 \\
& \leq & \frac{[w]_{RH_2^d}^4}{|J|} \sum_{I\in\mathcal{D}(J)} |I| \big (\langle w\rangle_{I_+}^2 + \langle w\rangle_{I_-}^2\big ) \left |\frac{\Delta_Iw}{\langle w\rangle_I}\right |^2  +
\frac{1}{|J|}\sum_{I\in\mathcal{D}(J)} |I| \left | \frac{\Delta_I(w^2)}{\langle w\rangle_{I}}\right |^2.
\end{eqnarray*}
Now using that $a^2+b^2\leq (a+b)^2$ for $a,b\geq 0$, and that $\langle w\rangle_{I_+}+\langle w\rangle_{I_-} = 2\langle w\rangle_I$ on the first sum, and    Proposition~\ref{lem:CHMPW}, with $p=2$, to estimate the second sum, we get
\begin{eqnarray*}
\Sigma_1 & \lesssim & \frac{[w]_{RH_2^d}^4}{|J|} \sum_{I\in\mathcal{D}(J)} |I| \langle w\rangle_{I}^2 \left |\frac{\Delta_Iw}{\langle w\rangle_I}\right |^2
+ C\langle w^2\rangle_J\\
& \leq & C[w]_{RH_2^d}^4 \langle w\rangle_J^2 + C\langle w^2\rangle_J \leq C_w \langle w^2\rangle_J.
\end{eqnarray*}
Where in the last line we used first, that $w\in RH_2^d$ and Buckley's characterization of $RH_2^d$ (Proposition~\ref{prop:Buckley}),  and second, H\"older's inequality. 
\end{proof}
The constant $C_w$ we obtain is of order $[w]_{RH_2^d}^6$ as the constant in Buckley's characterization is comparable to $[w]_{RH_2^d}^2$, and the constant in Lemma~\ref{lem:CHMPW}, when $w\in RH_2^d$, is comparable to $[w]_{RH_2^d}^2[w^2]_{RH_1}$. In other words, this argument yields a bound of the type $\|S_w\|_{L^2\to L^2}\leq C [w]_{RH_2^d}^{3}$ which is far from the optimal quadratic bound \cite{Per09}.


Notice that this argument does not require $w$ to be dyadic doubling or $w\in A_{\infty}^d$.
Corollary~\ref{coro:Per}, in the case that $w$ is a dyadic doubling weight, is a particular case of the more general one-weight result \cite[Theorem 7, proof in Section 3.3]{Per94}, that we state now.

\begin{thm}[\cite{Per94}]\label{thm:Per} Given $\sigma \in A_{\infty}^d$, and  a dyadic $\sigma$-doubling weight\footnote{A weight $\omega$ is  dyadic $\sigma$-doubling  if  $D^{\sigma}(w):=\sup_{I\in\mathcal{D}}\frac{\int_{\widetilde{I}} \omega\,d\sigma}{\int_I \omega \, d\sigma} <\infty$. Note that necessarily $D^{\sigma}(w)\geq 2$.} $\omega$. Then  $\omega \in A_p^d(d\sigma )$ if and only if $S^d_{\sigma}$ is bounded on $L^p(\omega d\sigma)$.
\end{thm}
The proof follows Buckley's proof in \cite{Buc93} when $d\sigma=dx$, for $p=2$, and then runs an extrapolation argument to get $1<p<\infty$. The linear bound when $p=2$, under the weaker assumption that $\sigma$ is dyadic doubling (not necessarily $\sigma\in A_{\infty}^d$), can be found in \cite[Theorem 4]{Per09}. More precisely, it is shown that $\|S_{\sigma}^df\|_{L^2(\omega d\sigma)}\leq C[\omega]_{A_2^d(d\sigma)}\|f\|_{L^2(\omega d\sigma)}$, where the constant $C$ depends on the dyadic doubling constant of $\sigma$, the proof uses a Bellman function argument.

Corollary~\ref{coro:Per},   in the case that $w$ is a dyadic doubling weight, is obtained from Theorem~\ref{thm:Per} by choosing $d\sigma=wdx$ and $\omega=w^{-1}$ so that $\omega d\sigma=dx$,  recalling that $w^{-1}\in A^d_p(w)$ if and only if $w\in RH_{p'}^d$, $\frac1p+\frac{1}{p'}=1$ (see Remark~\ref{remark1}),  and setting $p=2$. If we choose   $d\sigma=wdx$ and $\omega = uw^{-1}$ in Theorem~\ref{thm:Per} (under the weaker assumption in \cite[Theorem 4]{Per09}), we will get the following one weight theorem for the dyadic square function.

\begin{coro}\label{coro:Per2} 
Given a dyadic doubling weight $w$, and a weight $u$, with $uw^{-1} \in A_2^d(w)$ then the weighted dyadic square function $S_w^d$ is bounded on $L^2(u)$.
\end{coro}
To deduce  this corollary from Theorem~\ref{thm:square-function-two-weight}, we will use the following weight lemma which appeared in \cite[Lemma 10]{Per09}, stated in a way useful for us.
\begin{lem}[\cite{Per09}]\label{lem:WeightLemma} If $u^{-1}w\in A_2(w)$, and $w$ is a dyadic doubling weight, then
\[ \frac{1}{|J|} \sum_{I\in \mathcal{D}(J)} \frac{|\langle u^{-1}w, h_I^w\rangle_w|^2}{\langle u^{-1}w^2\rangle_I} \langle w\rangle_I\leq 18 D^{3}(w) [u^{-1}w]_{A_2(w)} \langle u^{-1}w^2\rangle_J.\]
Where $D(w)$ is the dyadic doubling constant of $w$.
\end{lem}
This lemma was proven using a Bellman function argument.

\begin{proof}[Proof of Corollary~\ref{coro:Per2}]
    If $u=v$ then condition (i) in Theorem~\ref{thm:square-function-two-weight} is  $uw^{-1} \in A_2^d(w)$. 
  Condition (ii) in this case reads
  \begin{equation}\label{eq:uwinvA2(w)}
 \Sigma_2\simeq \frac{1}{|J|}\sum_{I\in\mathcal{D}(J)} |I| \big |\Delta_I^w(u^{-1}w)\big |^2
 \langle u\rangle_I \leq C { \langle u^{-1}w^2\rangle_J}. 
  \end{equation}
Where $\Sigma_2$ is the left-hand-side of (ii), and the similarity occurs because we are assuming $w$ is dyadic doubling.
  
  We will show that if $uw^{-1} \in A_2^d(w)$ then \eqref{eq:uwinvA2(w)} holds.  Note that $uw^{-1}\in A_2^d(w)$ if and only if $u^{-1}w\in A^d_2(w)$, moreover $[uw^{-1}]_{A_2^d( w)}=[u^{-1}w]_{A_2^d(w)}$.  We can now use Lemma~\ref{lem:WeightLemma} to estimate $\Sigma_2$. Indeed, first, using \eqref{eq:DeltaIw(f)}; second,  using that  $\langle u\rangle_I\leq [uw^{-1}]_{A^d_2(w)} \langle w\rangle_I^2 /\langle u^{-1}w^2\rangle_I$, since  $uw^{-1} \in A_2^d(w)$; third, using that $w$ is dyadic doubling; and fourth, using Lemma~\ref{lem:WeightLemma}, we get
  \begin{eqnarray*}
  \Sigma_2 & = &\frac{1}{|J|}\sum_{I\in\mathcal{D}(J)} |I| \frac{w(I)}{w(I_+)w(I_-)} |\langle u^{-1}w, h_I^w\rangle_w|^2 \langle u\rangle_I \\
  & \leq & [uw^{-1}]_{A^d_2(w)} \frac{1}{|J|}\sum_{I\in\mathcal{D}(J)} |I| \frac{w(I)}{w(I_+)w(I_-)} |\langle u^{-1}w, h_I^w\rangle_w|^2 \frac{\langle w\rangle_I^2}{\langle u^{-1}w^2\rangle_I}\\
  & \leq & D(w)^2 [uw^{-1}]_{A^d_2(w)} \frac{1}{|J|}\sum_{I\in\mathcal{D}(J)}   \frac{|\langle u^{-1}w, h_I^w\rangle_w|^2}{\langle u^{-1}w^2\rangle_I}\langle w\rangle_I \\
  &\leq & 18 D(w)^5[uw^{-1}]^2_{A^d_2(w)} \langle u^{-1}w^2\rangle_J.
  \end{eqnarray*}
  
\end{proof}
This argument yields the following linear bound on the $A_2$-characteristic times a power of the doubling constant of $w$,
\[ \|S_w^d(f)\|_{L^2(u)\to L^2(u)} \lesssim D(w)^{5/2} [uw^{-1}]_{A_2^d(w)}.\]

The state of the art one-weight estimate for the dyadic square function is a linear bound on the $A_2$-characteristic of the weight even in the non-homogeneous setting and when the dyadic filtration is not regular \cite{DIPTV19}. That means no doubling constant should appear. Note, that in higher dimensions, the doubling constant is a dimensional constant. Having estimates that do not depend on the doubling constant, produces results that could potentially be  transferred to infinite dimensional settings (eg: UMD spaces). In \cite{DIPTV19}  they also show that the square root in the $A_2$-characteristic lower bound for the dyadic square function
\[ \|S^d\|_{L^2(w)\to L^2(w)} \geq C \sqrt{[w]_{A_2^d}},\]
known in the homogeneous case with a regular filtration  \cite{PP02, Per09}, no longer holds in this more general setting. Instead a linear bound on the $A_2$-characteristic holds. It is interesting that in the probability community, weighted inequalities for the martingale differences square function were known when the weight was in the corresponding $A_2$ class back in 1979  \cite{BL79}.


%
%

\section{Dyadic square function meets Haar multiplier}\label{sec:tHaar}

Recently in \cite[Theorem 1]{CHMPW24}, necessary and sufficient conditions for the two-weight boundedness of $T^t_{w,\sigma}$, the signed $t$-Haar multipliers (defined in \eqref{eq:t-HaarMultiplier}), were found. When restricting to $t=1$, this refers to the operator 
\[T_{w,\sigma}f(x)= \sum_{I\in\mathcal{D}} \sigma_I \frac{w(x)}{\langle w\rangle_I} \langle f,h_I\rangle h_I(x), \quad\mbox{where}\; \sigma_I=\pm 1.\]
We record the result in the case $t=1$, to guide the discussion.

\begin{thm}[Two-weight theorem \cite{CHMPW24}] \label{thm:two-weight} 
Given a triple of weights $(u,v,w)$. Denote $\Delta_I w:= \langle w\rangle_{I_+} - \langle w\rangle_{I_-} $. The  signed Haar multipliers $T_{w,\sigma}$ are uniformly (in $\sigma$) bounded from $L^2(u)$ into $L^2(v)$ if and only if the following four conditions hold,
\begin{itemize}
\item[{\rm (i)}] Joint three-weight condition: $\displaystyle{\quad C_1:=\sup_{I\in\mathcal{D}} \frac{\langle u^{-1}\rangle_I \langle v w^{2}\rangle_I}{\langle w\rangle_I^{2}} <\infty}$. 
\item[{\rm (ii)}] Carleson condition: there is a constant $C_2>0$ such that for all $J\in\mathcal{D}$, 
$$\frac{1}{|J|} \sum_{I\in\mathcal{D}(J)}  |I||\Delta_Iu^{-1}|^2 \frac{\langle vw^{2}\rangle_I}{\langle w\rangle_I^{2}}\leq C_2\,\langle u^{-1}\rangle_J.$$
\item[ {\rm (iii)}] Dual Carleson condition: there is a  constant $C_3>0$ such that for all $J\in\mathcal{D}$, 
$$\frac{1}{|J|} \sum_{I\in\mathcal{D}(J)}  |I||\Delta_I (vw^{2})|^2 \frac{\langle u^{-1}\rangle_I}{\langle w\rangle_I^{2}}\leq C_3\, \langle vw^{2} \rangle_J.$$
\item[{\rm (iv)}] The positive operator $P_{w,\lambda}$ is bounded from $L^2(u)$ into $L^2(v)$ with operator norm $C_4>0$, where  
$$P_{w,\lambda}f(x) := \sum_{I\in\mathcal{D}} \frac{w(x)}{|I|}\lambda_I \langle f\rangle_I \mathbbm{1}_I(x)\,\, \text{and} \,\, \lambda_I= \frac{|\Delta_I(u^{-1})|}{\langle u^{-1}\rangle_I}\frac{|\Delta_I(vw^2)|}{\langle vw^2 \rangle_I} \frac{|I|}{\langle w\rangle_I}. $$
\end{itemize}
  Moreover  $ \|T_{w,\sigma} \|_{L^2(u)\to L^2(v)}\lesssim  \sqrt{C_1} + \sqrt{C_2} +\sqrt{C_3}+ C_4.$
\end{thm}

In \cite{Per94}, it was proven that $T_w$ being bounded on $L^p(\R )$ is equivalent to $w \in RH_p^d$ and also equivalent to $S_w$ being bounded on $L^{p'}(\R )$ with \( \frac{1}{p}+\frac{1}{p'}=1 \). After proving  Theorem~\ref{thm:two-weight}, it is a natural question to see if some version of this equivalence holds in the two-weight context. We will impose some conditions on the weights in order to pursue that, let us define the quantities:
\begin{align}
Q_J(u,v,w):=& \frac{1}{|J|}\sum_{I\in \mathcal{D}(J):I\neq J} |\widetilde{I}| \left | \frac{\langle u^{-1}\rangle_{I}}{\langle w\rangle_{I}}\right |^2\left | 1-\frac{\langle w\rangle_{I}}{\langle w\rangle_{\widetilde{I}}}\right |^2 \langle vw^2\rangle_{\widetilde{I}}. \label{def:Q} \\
R_J(u,v,w):=& \frac{1}{|J|}\sum_{I\in \mathcal{D}(J):I\neq J} |\widetilde{I}| \left | \frac{\langle vw^2 \rangle_{I}}{\langle w\rangle_{I}}\right |^2\left | 1-\frac{\langle w\rangle_{I}}{\langle w\rangle_{\widetilde{I}}}\right |^2 \langle u^{-1}\rangle_{\widetilde{I}}.  \label{def:R}
\end{align}

\begin{thm}\label{t-HaarMultipliers-twoweight}
 Given a triple of weights $(u,v,w)$, assume that  $w$ is a dyadic doubling weight, 
 $Q_J(u,v,w) \leq  C_1\langle u^{-1}\rangle_{J}$ and $R_J(u,v,w) \leq  C_2\langle vw^2\rangle_{J}$ , for all $J \in \mathcal{D}$.   Then the Haar Multipliers $T_{w,\sigma}$ are uniformly bounded (on $\sigma$) from $L^2(u)$ to $L^2(v)$ if and only if the following three conditions hold simultaneously:
    \begin{itemize}
        \item[(i)] $S^d_w$ is bounded from $L^2(uw^2)$ to $L^2(vw^2)$, 
        \item[(ii)] $S^d_w$ is bounded from $L^2(v^{-1})$ into $L^2(u^{-1})$, 
        \item[(iii)] The positive operator $P_{w,\lambda}$ is bounded from $L^2(u)$ into $L^2(v)$, where
        $$ P_{w,\lambda}f(x):= \sum_{I \in \mathcal{D}}  \frac{w(x)}{|I|} \lambda_I \langle f \rangle_I \chi_I(x),  \,\, \text{and} \,\, \lambda_I= \frac{|\Delta_I(u^{-1})|}{\langle u^{-1}\rangle_I}\frac{|\Delta_I(vw^2)|}{\langle vw^2 \rangle_I} \frac{|I|}{\langle w\rangle_I}. $$
        \end{itemize}
\end{thm}

\begin{remark}
In the case $w=1$ the operator becomes $T_{\sigma}$,  the martingale transform. Moreover, both $Q_J(u,v,1)=R_J(u,v,1)=0$ and then the hypothesis of the above Theorem are trivially satisfied (namely, $w\in A_{\infty}^d$, $Q_J(u,v,w) \leq  C_1\langle u^{-1}\rangle_{J}$ and $R_J(u,v,w) \leq  C_2\langle vw^2\rangle_{J}$, for all $J \in \mathcal{D}$). In this case, the Theorem will exactly match the result \cite[Corollary 3.9]{BCMP17}.
\end{remark}
 
\begin{remark} In the case $u=v=1$ and $w$ a dyadic doubling weight,  we recover the equivalence of $T_w$, $S_w$ and $w \in RH^d_2$ presented in \cite{Per94}. Indeed, when $w\in RH_2^d$ and dyadic doubling, then it is in $A_{\infty}^d$, and  the conditions on $Q_J$ and $R_J$ hold,
\begin{eqnarray*} Q_J(1,1,w) & \lesssim & \frac{1}{|J|}\sum_{I\in \mathcal{D}(J):I\neq J} |\widetilde{I}| \left | \frac{\Delta_{\widetilde{I}} w}{\langle w\rangle_{\widetilde{I}}}\right |^2 
 \frac{\langle w^2\rangle_{\widetilde{I}}}{\langle w\rangle_{\widetilde{I}}^2} \, \leq \, C [w]_{RH_2^d}^2\log [w]_{A_{\infty}^d} =C_w,\\
R_J(1,1,w) & = & \frac{1}{|J|}\sum_{I\in \mathcal{D}(J):I\neq J} |\widetilde{I}| \left | \frac{\Delta_{\widetilde{I}} w}{\langle w\rangle_{\widetilde{I}}}\right |^2
 \frac{\langle w^2 \rangle^2_{I}}{\langle w\rangle^4_{I}} \langle w\rangle_{\widetilde I}^2
\, \leq \, C[w]_{RH_2}^4  \langle w\rangle^2 =C_w \langle w^2\rangle_J.
\end{eqnarray*}
In the estimate for  $Q_J$, we used Proposition~\ref{thm:FKP}, and for $R_J$, we used Proposition~\ref{prop:Buckley} \textnormal{(ii)}.
\end{remark}


\begin{proof}[Proof of Theorem~\ref{t-HaarMultipliers-twoweight}]
We want to show that the conditions in Theorem~\ref{t-HaarMultipliers-twoweight} imply the necessary and sufficient conditions in Theorem~\ref{thm:two-weight}, and vice versa.

First, note  that condition (iii) in Theorem~\ref{t-HaarMultipliers-twoweight} is identical to condition (iv) in Theorem~\ref{thm:two-weight}. We have to focus on conditions (i) and (ii) in Theorem~\ref{t-HaarMultipliers-twoweight} being equivalent to conditions (i), (ii), and (iii) in Theorem~\ref{thm:two-weight}.

Second, since we assumed $w \in A^d_{\infty}$, then Corollary~\ref{coro:iff}  provides the following  necessary and sufficient conditions for the boundedness of $S_w$ from $L^2(uw^2)$ into $L^2(vw^2)$:
\begin{itemize}
    \item[(a)] $(uw,vw) \in A^d_2(w)$,
    \item[(b)] $ \displaystyle \frac{1}{w(J)} \sum_{I \in \mathcal{D}(J)} \big | \Delta^w_I (u^{-1}w^{-1}) \big|^2 \left [ \frac{w(I^-)^2(vw^2)(I^+) + w(I^+)^2(vw^2)(I^-)}{w(I)^2} \right ]  \leq C \langle u^{-1}w^{-1}\rangle_J^w$,
\end{itemize}

Condition (a) is exactly  condition (i) in Theorem~\ref{thm:two-weight}, by Lemma~\ref{prop:equiv-3weight-conditions}.

Since $w$ is a dyadic doubling weight, condition (b) can be rewritten in terms of Lebesgue averages as
\begin{equation}\label{condition-(i)-uw2-vw2}
\frac{1}{|J|} \sum_{I \in \mathcal{D}(J) } |I|  \bigg|\frac{\langle u^{-1}\rangle_{I^+}}{\langle w\rangle_{I^+}} -\frac{\langle u^{-1}\rangle_{I^-}}{\langle w\rangle_{I^-}}\bigg|^2   \langle vw^2\rangle_I\leq C \langle u^{-1}\rangle_J.
\end{equation}
This is very  close to  condition (ii) in Theorem~\ref{thm:two-weight}, namely
\begin{equation}\label{condition(ii)t-Haar}
\frac{1}{|J|} \sum_{I \in \mathcal{D}(J) } |I|  \bigg|\frac{\langle u^{-1}\rangle_{I^+}}{\langle w\rangle_{I}} -\frac{\langle u^{-1}\rangle_{I^-}}{\langle w\rangle_{I}}\bigg|^2 \langle vw^2\rangle_I  \leq C \langle u^{-1}\rangle_J.
\end{equation}
 Let us denote by
 \begin{align}
 \Delta_{I,w}(u^{-1})&:=\frac{\langle u^{-1}\rangle_{I^+}}{\langle w\rangle_{I^+}} -\frac{\langle u^{-1}\rangle_{I^-}}{\langle w\rangle_{I^-}}, \label{eq:delta_original} \\
 \Delta_{I,w}'(u^{-1})&:=\frac{\langle u^{-1}\rangle_{I^+}}{\langle w\rangle_{I}} -\frac{\langle u^{-1}\rangle_{I^-}}{\langle w\rangle_I}. \label{eq:delta_prime}
 \end{align}
  
With this notation,   \eqref{condition-(i)-uw2-vw2}  and \eqref{condition(ii)t-Haar} become,
\begin{equation}
\frac{1}{|J|} \sum_{I \in \mathcal{D}(J) } |I| \, |\Delta_{I,w}(u^{-1})|^2 \langle vw^2\rangle_I  \leq C \langle u^{-1}\rangle_J. \label{condition-(i)-uw2-vw2_modified}
\end{equation}
\begin{equation}
\frac{1}{|J|} \sum_{I \in \mathcal{D}(J) } |I| \, |\Delta_{I,w}'(u^{-1})|^2 \langle vw^2\rangle_I  \leq C \langle u^{-1}\rangle_J, \label{condition(ii)t-Haar_modified}
\end{equation}

Assume  \eqref{condition(ii)t-Haar_modified} and let us show \eqref{condition-(i)-uw2-vw2_modified}. Starting from \eqref{eq:delta_original} and adding and subtracting $\langle u^{-1}\rangle_{I^+}/\langle w\rangle_I$ and $\langle u^{-1}\rangle_{I^-}/\langle w\rangle_I$, and observing that $|a+b+c|^2\leq 3(|a|^2+|b|^2+|c|^2)$, we get
\begin{eqnarray*}
 |\Delta_{I,w}(u^{-1})|^2&\lesssim&  
 \bigg|\frac{\langle u^{-1}\rangle_{I^+}}{\langle w\rangle_{I^+}} -\frac{\langle u^{-1}\rangle_{I^+}}{\langle w\rangle_{I}}\bigg|^2 +  |\Delta_{I,w}'(u^{-1})|^2+ \bigg |\frac{\langle u^{-1}\rangle_{I^-}}{\langle w\rangle_{I}} -\frac{\langle u^{-1}\rangle _{I^-}}{\langle w\rangle _{I^-}}\bigg|^2\\
 &=&  |\Delta_{I,w}'(u^{-1})|^2+  \left | \frac{\langle u^{-1}\rangle_{I^+}}{\langle w\rangle_{I^+}}\right |^2\left | 1-\frac{\langle w\rangle_{I^+}}{\langle w\rangle_I}\right |^2+ \left | \frac{\langle u^{-1}\rangle_{I^-}}{\langle w\rangle_{I^-}}\right |^2\left | 1-\frac{\langle w\rangle_{I^-}}{\langle w\rangle_I}\right |^2.
 \end{eqnarray*}
 Adding over $I\in \mathcal{D}(J)$ , we get
\begin{eqnarray*}
 \frac{1}{|J|} \sum_{I \in \mathcal{D}(J) } |I|\, |\Delta_{I,w}(u^{-1})|^2 \langle vw^2\rangle_I  \leq \; \frac{1}{|J|} \sum_{I \in \mathcal{D}(J) } |I| \, |\Delta_{I,w}'(u^{-1})|^2 \langle vw^2\rangle_I + Q_J(u,v,w).
 \end{eqnarray*}
 Now, using \eqref{condition(ii)t-Haar_modified} and that $ Q_J(u,v,w) \leq C \langle u^{-1} \rangle_J$ we have that \eqref{condition-(i)-uw2-vw2_modified} is proven:
 \begin{eqnarray*}
 \frac{1}{|J|} \sum_{I \in \mathcal{D}(J) } |I|\, |\Delta_{I,w}(u^{-1})|^2 \langle vw^2\rangle_I  \leq  \; C \langle u^{-1} \rangle_J.
 \end{eqnarray*}
 
 Likewise,  now assume \eqref{condition-(i)-uw2-vw2_modified} and let us show \eqref{condition(ii)t-Haar_modified}. 
 Adding and subtracting $\langle u^{-1}\rangle_{I^+}/\langle w\rangle_{I^+}$ and $\langle u^{-1}\rangle_{I^-}/\langle w\rangle_{I^-}$, and observing that $|a+b+c|^2\leq 3(|a|^2+|b|^2+|c|^2)$, we get
\begin{eqnarray*}
 |\Delta'_{I,w}(u^{-1})|^2&\lesssim&  
  |\Delta_{I,w}(u^{-1})|^2+  \left | \frac{\langle u^{-1}\rangle_{I^+}}{\langle w\rangle_{I^+}}\right |^2\left | 1-\frac{\langle w\rangle_{I^+}}{\langle w\rangle_I}\right |^2+ \left | \frac{\langle u^{-1}\rangle_{I^-}}{\langle w\rangle_{I^-}}\right |^2\left | 1-\frac{\langle w\rangle_{I^-}}{\langle w\rangle_I}\right |^2.
 \end{eqnarray*}
 With the same notation as before, adding over $I\in \mathcal{D}(J)$, we get
\begin{eqnarray*}
 \frac{1}{|J|} \sum_{I \in \mathcal{D}(J) } |I|\, |\Delta'_{I,w}(u^{-1})|^2 m_I (vw^2)  \leq \; \frac{1}{|J|} \sum_{I \in \mathcal{D}(J) } |I| \, |\Delta_{I,w}(u^{-1})|^2 m_I (vw^2) + Q_J(u,v,w).
 \end{eqnarray*}
Again, using \eqref{condition-(i)-uw2-vw2_modified} and the hypothesis over $Q_J(u,v,w)$, we obtain \eqref{condition(ii)t-Haar_modified}:
 \begin{eqnarray*}
 \frac{1}{|J|} \sum_{I \in \mathcal{D}(J) } |I|\, |\Delta'_{I,w}(u^{-1})|^2 \langle vw^2\rangle_I  \leq \; C \langle u^{-1}\rangle_J.
 \end{eqnarray*}

Third, since we assumed $w \in A^d_{\infty}$, then Corollary~\ref{coro:iff}  provides  the following  necessary and sufficient conditions for the boundedness of $S_w$ from $L^2(v^{-1})$ into $L^2(u^{-1})$:

\begin{itemize}
    \item[(c)] $(v^{-1}w^{-1},u^{-1}w^{-1}) \in A^d_2(w)$,
    \item[{\rm (d)}]  $\; \displaystyle \frac{1}{w(J)} \sum_{I \in \mathcal{D}(J) } \big| \Delta^w_I (vw) \big|^2  \left [ \frac{w(I^-)^2u^{-1}(I^+) + w(I^+)^2u^{-1}(I^-)}{w(I)^2} \right ]  \leq C \langle vw\rangle_J^w$.
\end{itemize}

Condition (c) is equivalent to condition (i) in Theorem~\ref{thm:two-weight} by Proposition~\ref{prop:equiv-3weight-conditions}.
When $w$ is dyadic doubling, condition (d) is equivalent to

\hskip .1in (d') $\quad \displaystyle \frac{1}{w(J)} \sum_{I \in \mathcal{D}(J) }w(I) \big| \Delta^w_I (vw) \big|^2 \langle u^{-1}w^{-1}\rangle^w_I  \leq C \langle vw\rangle_J^w$.\\
Condition (c') can be written in terms of Lebesgue averages as,
\begin{equation}\label{condition(b)dual-spaces}
\frac{1}{|J|} \sum_{I \in \mathcal{D}(J) } |I| \bigg|\frac{\langle vw^2\rangle_{I^+}}{\langle w\rangle_{I^+}} -\frac{\langle vw^2\rangle_{I^-}}{\langle w\rangle_{I^-}}\bigg|^2 \langle u^{-1}\rangle_I\leq C \langle vw^2\rangle_J.
\end{equation}
This is almost condition (iii) in Theorem~\ref{thm:two-weight}, namely
\begin{equation}\label{condition(iii)t-Haar}
\frac{1}{|J|} \sum_{I \in \mathcal{D}(J) } |I| \bigg|\frac{\langle vw^2\rangle_{I^+}}{\langle w\rangle_{I}} -\frac{\langle vw^2\rangle_{I^-}}{\langle w\rangle_{I}}\bigg|^2 \langle u^{-1}\rangle_I\leq C \langle vw^2\rangle_J.
\end{equation}

Analogously as we proved that \eqref{condition-(i)-uw2-vw2} and \eqref{condition(ii)t-Haar} are equivalent when $Q_J(u,v,w) < C \langle u^{-1}\rangle_J \,$ for all  $\, J \in \mathcal{D}$, we can prove that \eqref{condition(b)dual-spaces} and \eqref{condition(iii)t-Haar} are equivalent assuming  $R_J(u,v,w) \leq  C \langle vw^2\rangle_J\,$  for all  $\, J \in \mathcal{D}$.

 \end{proof}

Alternately, we could show that if there is a constant $C>0$ such that 
\begin{align}
 \frac{1}{|J|} & \sum_{I\in \mathcal{D}(J)} |I| \left |\frac{\Delta_Iw}{\langle w\rangle_I}\right |^2 \langle u^{-1}\rangle_I \leq C \langle u^{-1}\rangle_J \quad \forall J \in \mathcal{D} \label{eq:u(-1)-Carleson}, \\
 \frac{1}{|J|} & \sum_{I\in \mathcal{D}(J)} |I| \left |\frac{\Delta_Iw}{\langle w\rangle_I}\right |^2 \langle vw^2\rangle_I \leq C \langle vw^2\rangle_J  \quad  \forall J \in \mathcal{D}, \label{eq:vw2-Carleson} 
 \end{align}
 then $Q_J(u,v,w) \lesssim \langle u^{-1}\rangle_{J}$ and $R_J(u,v,w) \lesssim \langle vw^2\rangle_{J}$ hold  for all $J \in \mathcal{D}$.
 
  \begin{lem}\label{lem:QJuvw}  Assume that $w$ is dyadic doubling, that $(uw,vw)\in A_2^d(w)$, and \eqref{eq:u(-1)-Carleson}, then $Q_J(u,v,w)\lesssim \langle u^{-1}\rangle_J$.
 \end{lem}

 \begin{proof}[Proof of Lemma~\ref{lem:QJuvw}]
 First, note that $|\Delta_{\widetilde{I}} w|= \frac12 |\langle w\rangle_{\widetilde{I}}-\langle w\rangle_I|$, hence
 \begin{eqnarray*}
 Q_J(u,v,w) & = & \frac{2}{|J|}\sum_{I\in \mathcal{D}(J):I\neq J} |\widetilde{I}| \left | \frac{\langle u^{-1}\rangle_{I}}{\langle w\rangle_{I}}\right |^2 \left | \frac{\Delta_{\widetilde{I}} w}{\langle w\rangle_{\widetilde{I}}}\right |^2 \langle vw^2\rangle_{\widetilde{I}} \\
 & \lesssim &  \frac{1}{|J|}\sum_{I\in \mathcal{D}(J):I\neq J} |\widetilde{I}|  \frac{\langle u^{-1}\rangle_{\widetilde{I}}\langle vw^2\rangle_{\widetilde{I}}}{\langle w\rangle^2_{\widetilde{I}}} \left | \frac{\Delta_{\widetilde{I}} w}{\langle w\rangle_{\widetilde{I}}}\right |^2  \langle u^{-1}\rangle_{\widetilde{I}}\\
 &\lesssim & [uw,vw]_{A_2(w)}  \frac{1}{|J|}\sum_{I\in \mathcal{D}(J):I\neq J} |\widetilde{I}|  \left | \frac{\Delta_{\widetilde{I}} w}{\langle w\rangle_{\widetilde{I}}}\right |^2  \langle u^{-1}\rangle_{\widetilde{I}}\\
 & \lesssim & [uw,vw]_{A_2(w)} \langle u^{-1}\rangle_{J}.
 \end{eqnarray*}
 Where we used in the first inequality the facts  that  $w$ is dyadic doubling, and that $\langle u^{-1}\rangle_I \leq 2\langle u^{-1}\rangle_{\widetilde{I}}$. We used the assumption that $(uw,vw)\in A_2^d(w)$ and  \eqref{eq:u(-1)-Carleson} in the last inequality.  
 \end{proof}
 
  \begin{lem}\label{lem:RJuvw}  Assume that $w$ is dyadic doubling, that $(uw,vw)\in A_2^d(w)$, and \eqref{eq:vw2-Carleson}, then $R_J(u,v,w)\lesssim \langle vw^2\rangle_J$.
 \end{lem}
 The proof is analogous to the proof of Lemma~\ref{lem:QJuvw}. Therefore, in Theorem~\ref{t-HaarMultipliers-twoweight} we can replace the three-weight hypotheses, 
$Q_J(u,v,w) \leq  C_1\langle u^{-1}\rangle_{J}$ and $R_J(u,v,w) \leq  C_2\langle vw^2\rangle_{J}$, for all $J \in \mathcal{D}$, by  the hypotheses~\eqref{eq:u(-1)-Carleson} and~\eqref{eq:vw2-Carleson}, where the roles of $u$ and $v$ have been decoupled, and the theorem still holds true.  This decoupling indicates that this version of the theorem could be weaker than Theorem~~\ref{t-HaarMultipliers-twoweight}.

Conditions \eqref{eq:u(-1)-Carleson} and \eqref{eq:vw2-Carleson} are weighted versions of the celebrated R. Fefferman, C. Kenig, J. Pipher characterization of dyadic $A_{\infty}$ \cite{FKP91},  stated in Section~\ref{sec:prelim} as Proposition~\ref{thm:FKP}.

Proposition~\ref{thm:FKP}  has far reaching consequences and generalization described eloquently in another paper in this collection \cite{Pet25}. 

 In the language of Carleson sequences, Proposition~\ref{thm:FKP} states that if $w\in A_{\infty}^d$ then the sequence $\{b_I^2\}_{I\in\mathcal{D}}$, where $b_I^2=  |I| \left |\frac{\Delta_Iw}{\langle w\rangle_I}\right |^2$, is a Carleson sequence, or equivalently, the function $b=\sum_{I\in\mathcal{D}} b_I h_I$ is in dyadic $BMO^d$.
 Condition \eqref{eq:u(-1)-Carleson} states that, in addition, the sequence $\{b_I^2\}_{I\in\mathcal{D}}$ is a $u^{-1}$-Carleson sequence, as defined in \cite{Per94}. We can deduce from \cite[Lemma 4]{Per94} the following lemma.
 \begin{lem}
 Assume that  $w \in A^d_{\infty}$. Then
 \begin{itemize}
 \item[(i)] if $ u^{-1} \in A^d_{\infty}(w)$ then \eqref{eq:u(-1)-Carleson}  holds.
 \item[(ii)] if $ vw^2 \in A^d_{\infty}(w)$ then \eqref{eq:vw2-Carleson} holds.
 \end{itemize}
 \end{lem}
Here $u^{-1}\in A_{\infty}^d(w)$ means that the sequence $\{a_I\}_{I\in\mathcal{D}}$, where $a_I=|I| \left | \frac{\langle u^{-1}\rangle_I^w - \langle u^{-1}\rangle_{\widetilde{I}}^w}{\langle u^{-1}\rangle_{\widetilde{I}}^w}\right |^2$, is a $w$-Carleson sequence, that is, $ \frac{1}{|J|} \sum_{I\in\mathcal{D}(J)} a_I \langle w\rangle_I \lesssim \langle w\rangle_J$ (see Section~\ref{sec:Carleson}).

Therefore, in Theorem~\ref{t-HaarMultipliers-twoweight} we can replace the three-weight hypotheses, 
$Q_J(u,v,w) \lesssim  \langle u^{-1}\rangle_{J}$ and $R_J(u,v,w) \lesssim \langle vw^2\rangle_{J}$ , for all $J \in \mathcal{D}$, by  the hypotheses $ u^{-1} \in A^d_{\infty}(w)$ and $ vw^2 \in A^d_{\infty}(w)$, and the theorem will still hold true. Note that, in this version of the theorem,  we are still assuming that $w\in A_{\infty}^d$, and the roles of $u$ and $v$ have been decoupled in the preliminary hypotheses, an indication that this could  be a weaker version of the theorem.

To finish we would like to highlight another connection to \cite{FKP91}, via the resolvent of the dyadic paraproduct \cite{Per94}, mentioned in the introduction.  Considering two-weight analogues for the resolvent of the paraproduct led us  to search for necessary and sufficient conditions relating  the two weight boundedness of the Haar multipliers $T_{w,\sigma}$ and the weighted dyadic square function $S_w$, and to do that we had to first understand better the square function, which are the results discussed in this paper.

Attempting to characterize Sobolev spaces on Lipschitz curves via some geometric second difference, in~\cite{Per96} the third author encountered an operator, a \emph{variable Haar multiplier},  that looked like
\[P_bf(x)  =  \sum_{I\in\mathcal{D}(I_0)} \prod_{J\in\mathcal{D}:J\subsetneq I} \big (1+b_J h_J(x)\big ) \, \langle f, h_I\rangle h_I(x).\]
Where the sequence $\{b_I\}_{I\in\mathcal{D}}$ was dictated by the geometry of the curve, and so were the dyadic intervals, and the Haar functions (weighted with respect to a measure given by the Lipschitz map).
Formally, unfolding the products in the sum, one gets a Neumann series
\[P_bf(x)  =  \sum_{n=0}^{\infty} \pi_b^nf(x) \quad\quad\quad \mbox{``paraseries'',}\]
 where  $\; \pi_bf(x):= \sum_{J\in\mathcal{D}(I_0)} \langle f\rangle_J b_Jh_J(x)$ is the \emph{dyadic paraproduct}.
If  $\pi_b$ is a contraction in $L^p(I_0)$  then
$P_bf(x)  = (I-\pi_b )^{-1}(f)$ is bounded on $L^p(I_0)$. But being a contraction imposed a
heavy constraint on $b$, not warranted by  geometry of the problem. A different argument was needed to get the boundedness of $P_b$, without requiring the paraproduct to be a contraction. 

Instead,  a remarkable product formula introduced in \cite[Section 3.18]{FKP91} was used to
rewrite $P_b$, in terms of  $w=\prod_{J\in \mathcal{D}(I_0)} \big (1+b_Jh_J\big )$, as
\[P_bf(x)  =  P_wf(x)= \sum_{I\in\mathcal{D}(I_0)} \frac{w(x)}{\langle w\rangle_I \big (1+b_Ih_I(x)\big )}\, \langle f, h_I\rangle h_I(x) ,\]
where $\langle w\rangle_I=\frac{1}{|I|}\int_I w(x)\, dx$, and $b_I=\langle b,h_I\rangle={\langle w, h_I\rangle}/{\langle w\rangle_I}$.  In this light, the operator $P_w$ seemed remarkably close to the operator $T_w$ discussed in this section. Indeed, in \cite{Per94}, it was shown that given a dyadic doubling weight $w$,  $P_w$ is bounded on $L^p$ if, and only if, $T_w$ is bounded on $L^p$ if, and only if, $S_w$ is bounded on $L^{p'}$ for $\frac1p+\frac{1}{p'} =1$, if, and only if $w\in RH_p^d$. This provided an alternative proof for the boundedness of the resolvent of the paraproduct, $P_b =(I-\pi_b)^{-1}$, unmasking hidden cancellations in the paraseries, and removing the geometric constraints imposed when assuming that the dyadic paraproduct was a contraction. 

 In \cite{FKP91}, necessary and sufficient   conditions on $b$ where found to ensure the convergence of the product to a doubling measure, that we record here.
  
 \begin{lem}[\cite{FKP91}, Lemma 3.20]\label{lem:FKP}
 Let $\{b_I\}_{I\in \mathcal{D(I_0)}}$ be a sequence satisfying ${|b_I |}/{\sqrt{|I|}}\leq 1-\epsilon$ for some $\epsilon>0$ and for all $I\in\mathcal{D}$. Then 
 the product $\prod_{I\in\mathcal{D}(I_0)} \big (1+b_Ih_I(x)\big )$ converges in  the weak* topology of measures to a  dyadic doubling probability measure $\mu$ on $I_0$. 
 
 Conversely, given a positive dyadic doubling probability measure $\mu$ on $I_0$, then the sequence $\{b_I\}_{I\in\mathcal{D}(I_0)}$ given by $b_I=\langle \mu,h_I\rangle/\langle \mu\rangle_I$ satisfies  $|b_I |/\sqrt{|I|}\leq 1-\epsilon $ for some $\epsilon>0$ and for all $I\in\mathcal{D}(I_0)$, furthermore the product $\prod_{I\in\mathcal{D}(I_0)} \big (1+b_Ih_I(x)\big )$ converges weak* to $\mu$.
  \end{lem}
  
   Proposition~\ref{thm:FKP} (\cite[Proposition 1]{FKP91}), can be paraphrased as ``$w$ is in $A_{\infty}^d$ if, and only if, $b \in BMO^d$". This and Lemma~\ref{lem:FKP},  established the first entries in a dictionary later extended by Buckley in \cite{Buc93} to $A_p^d$ and $RH_p^d$ (Proposition~\ref{prop:Buckley}), see also \cite{BR14}. For some applications to ``paraexponentials", see \cite{PW98}.
     This remarkable product formula has found applications in Data Science, providing the theoretical foundation for the representation of a data set as a measure in a very large hierarchically parametrized family of positive measures, whose parameters can be computed explicitly \cite{BJNS20}. 
     
     What we have described here are just three pages in the almost 60 pages long and very influential paper \cite{FKP91}. The paper includes a continuous version of these results with surprising  applications to elliptic PDEs. A hybrid (continuous-dyadic) version of the results described in these final paragraphs can be found \cite{GN02}.

\end{document}